\documentclass[a4paper,12pt]{amsart}

\usepackage[english]{babel}
\usepackage[latin1]{inputenc}
\usepackage{amsmath}
\usepackage{amsthm}
\usepackage{amsfonts}
\usepackage{amssymb, latexsym, graphics, showidx}
\usepackage{color}

\numberwithin{equation}{section}
\newtheorem{de}{Definition}[section]
\newtheorem{thm}{Theorem}[section]
\newtheorem{rem}[thm]{Remark}

\newtheorem{prop}[thm]{Proposition}
\newtheorem{lem}[thm]{Lemma}

\renewcommand{\dim}{\noindent\textbf{Proof.} }

\newcommand{\finedim}{{\unskip\nobreak\hfil\penalty50
   \hskip2em\hbox{}\nobreak\hfil\mbox{$\Box$ \qquad}
   \parfillskip=0pt \finalhyphendemerits=0\par\medskip}}
\newcommand{\R}{\mathbb{R}}

\newcommand{\Z}{\mathbb{Z}}
\newcommand{\Om}{\Omega}

\newcommand{\lam}{\lambda}
\newcommand{\al}{\alpha}

\newcommand{\xs}{\overline{x}}
\newcommand{\ys}{\overline{y}}
\newcommand{\ts}{\overline{t}}
\newcommand{\tas}{\overline{\tau}}
\newcommand{\zs}{\overline{z}}
\newcommand{\p}{\partial}
\newcommand{\s}{\sigma}
\newcommand{\ep}{\epsilon}
\newcommand{\I}{\mathcal{I}_1}

\newcommand{\Xs}{\overline{X}}
\newcommand{\Ys}{\overline{Y}}
\newcommand{\Zs}{\overline{Z}}
\newcommand{\os}{\overline{s}}
\newcommand{\yn}{y_{N+1}}
\newcommand{\zn}{z_{N+1}}

\newcommand{\beq}{\begin{equation}}
\newcommand{\eeq}{\end{equation}}
\newcommand{\beqs}{\begin{equation*}}
\newcommand{\eeqs}{\end{equation*}}
\newcommand{\beqa}{\begin{eqnarray}}
\newcommand{\eeqa}{\end{eqnarray}}
\newcommand{\beqas}{\begin{eqnarray*}}
\newcommand{\eeqas}{\end{eqnarray*}}

\title[]{Homogenization of the Peierls-Nabarro model\\ for dislocation dynamics}
\email{monneau@cermics.enpc.fr} \email{spatrizi@math.ist.utl.pt}

\begin{document}


\maketitle \centerline{\scshape R\'{e}gis Monneau }
\medskip
{\footnotesize
 \centerline{Universit\'e Paris-Est, CERMICS, Ecole des Ponts ParisTech,}
   \centerline{6-8 avenue Blaise Pascal, Cit\'{e}
Descartes, Champs sur Marne,}
   \centerline{77455 Marne la Vallée Cedex 2,
France}
} 

\medskip

\centerline{\scshape Stefania Patrizi}
\medskip
{\footnotesize
\centerline{Instituto Superior T\'{e}cnico,  Dep. de Matem\'{a}tica}
   \centerline{Av. Rovisco Pais Lisboa, Portugal}

}

\bigskip
\begin{abstract}This paper is concerned with a result of
homogenization of an integro-differential equation describing
dislocation dynamics. Our model involves both an anisotropic
L\'{e}vy operator of order 1 and a potential depending
periodically on $u/\ep$. The limit equation is a non-local
Hamilton-Jacobi equation, which is  an effective plastic law for
densities of dislocations moving in a single slip plane.
\end{abstract}

\section{Introduction}
In this paper we are interested in homogenization of the
Peierls-Nabarro model, which is a phase field model describing
dislocations. In this model a dislocation is described by a phase transition.
Dislocations are moving defects in crystals that can be described at several scales by different models:
\begin{itemize}
\item atomic scale (Frenkel-Kontorova model), 
\item microscopic scale (Peierls-Nabarro model), 
\item mesoscopic scale (Discrete dislocation dynamics), 
\item macroscopic scale (elasto-visco-plasticity with density of dislocations). 
\end{itemize}
Several changes of scales already exist in the literature: see for instance \cite{him} for a presentation of 
rigorous passages from atomic scale to microscopic scale, from microscopic scale to mesoscopic scale and from mescoscopic scale to macroscopic scale. Notice that the passage from Peierls-Nabarro model to the Discrete dislocation dynamics is only done in dimension 1 (see \cite{him} and \cite{gonzalezmonneau}). 
On the contrary in higher dimensions, the large scale limit of  a single 
phase transition described by the Peierls-Nabarro model shows that the line tension effect
is the much stronger term. The limit model appears to be 
the mean curvature motion (see \cite{is}). 

Our goal in this paper is to understand the large scale limit of the Peierls-Nabarro model in the case of a large number of phase transitions (i.e. of dislocations),
recovering at the limit a model with evolution of dislocation densities. In other words,  
we want to {\it perform a direct passage in any dimensions from the microscopic scale (Peierls-Nabarro model) to the 
macroscopic scale (elasto-visco-plasticity with density of dislocations)}.
In physics and mechanics, it is a great challenge to try to predict macroscopic elasto-visco-plasticity 
properties of materials (like metals), based on microscopic properties like dislocations.
In our work, we try to tackle this question in a very simplified geometry 
where all the dislocations are contained in the same slip plane with the same Burgers vector.
For a physical introduction to the Peierls-Nabarro
model, see for instance \cite{hl}; for a recent reference, see
\cite{WXM}; we also refer the reader to the paper of Nabarro
\cite{N} which presents an historical tour on the Peierls-Nabarro
model. See also Section \ref{physicalsec} for a more physical presentation of the Peierls-Nabarro model
and an interpretation of our results.

\subsection{Setting of the problem}

 The Peierls-Nabarro model has been originally introduced as a variational (stationary) model (see \cite{N}).
The time evolution Peierls-Nabarro model as a gradient flow dynamics has only been introduced 
quite recently, see for instance \cite{MBW} and \cite{Denoual}.
In the present paper we consider such a  time evolution  
Peierls-Nabarro model that can be written at the microscopic scale for the parameter $\ep=1$ as the following equation
\begin{equation}\label{uep}
\begin{cases}
\p_{t}
u^\epsilon=\I[u^\epsilon(t,\cdot)]-W'\left(\frac{u^\epsilon}{\epsilon}\right)+\s
\left(\frac{t}{\epsilon},\frac{x}{\epsilon}\right)&\text{in}\quad \R^+\times\R^N\\
u^\epsilon(0,x)=u_0(x)& \text{on}\quad \R^N.
\end{cases}
\end{equation}
For the physical application that we have in mind, we consider a three-dimensional crystal
which contains a crystallographic plane $\R^N$ with $N=2$.
This plane contains the dislocations that are represented by transitions of the phase function $u^\ep$.  
Here  $u^\ep$ solves the non local (and non linear) heat equation (\ref{uep}).
Indeed $\I$ stands here for an anisotropic half Laplacian (whose expression will be precised below).
Here the anisotropy comes both from the possible anisotropy of the elasticity of the crystal 
and from the fact that the Burgers vector is assumed to be contained in the slip plane $\R^N$ 
which creates a preferable direction.
The dynamics is assumed to be fully overdamped and 
then the right hand side of the equation 
is the sum of three force terms: $\I[u^\epsilon]$ is the elastic stress created by the dislocation themselves, $-W'$ is the force deriving from the potential $W$ describing the misfit between the two half crystals separated by the plane $\R^N$, and $\sigma$ is a stress created by the obstacles in the crystal or/and an applied exterior stress.
For simplicity $\sigma$ is assumed to be periodic 
in order to analyse by homogenization the effect on the dynamics of periodic obstacles everywhere in the crystal.
We consider time periodicity for two reasons:
one in order to take into account exterior periodic loads, and the second for generality.
Indeed, if $\sigma(t/\varepsilon,x/\varepsilon)$ is replaced by 
an oscillation at a different scale like $\sigma(t/\varepsilon^\gamma, x/\varepsilon^\gamma)$ with $\gamma\not= 1$,
then we expect (but it is not proven) that there is a two-scales homogenization effect.
If $\gamma>1$, then we expect that there is first homogenization of $\sigma$,
where only its mean value will be taken into account at the microscopic scale, 
and in a second step,
we get the macroscopic model by homogenization of the Peierls-Nabarro model with constant $\sigma$.
If $\gamma<1$, we expect first to freeze $\sigma$ and get the macroscopic  model  by homogenization of the Peierls-Nabarro model for constant $\sigma$, 
and in a second step we remind us that $\sigma$ is slowly oscillating,
and there is a second homogenization of the macroscopic model.

Here $\ep$ describes the ratio between the microscopic scale and the macroscopic scale,
and then is a small parameter. After a suitable rescaling at the macroscopic scale, 
the Peierls-Nabarro model becomes (\ref{uep}).
In this paper we investigate the limit as $\ep\rightarrow0$ of the viscosity
solution $u^\ep$ of (\ref{uep}).

We give the precise definitions and assumptions on the terms involved in (\ref{uep}).
Here $\I$ is an anisotropic L\'{e}vy operator of order 1, defined
on bounded $C^2$- functions for $r>0$ by
\beq\label{levy}\begin{split}\I[U](x)&=\int_{|z|\leq
r}(U(x+z)-U(x)-\nabla U(x)\cdot
z)\frac{1}{|z|^{N+1}}g\left(\frac{z}{|z|}\right)dz\\&
+\int_{|z|>r}(U(x+z)-U(x))\frac{1}{|z|^{N+1}}g\left(\frac{z}{|z|}\right)dz,\end{split}\eeq
where the function $g$ satisfies
\begin{itemize}
\item [(H1)] $g\in C({\bf S}^{N-1}),\,g> 0$,\,$g$ even.\end{itemize} On the functions $W$, $\sigma$ and $u_0$ we
assume:
\begin{itemize}
\item [(H2)]$W\in C^{1,1}(\R)$ and $W(v+1)=W(v)$ for any $v\in\R$;
\item [(H3)]$\sigma\in C^{0,1}(\R^+\times\R^N)$ and $\s(t+1,x)=\s(t,x)$, $\sigma(t,x+k)=\sigma(t,x)$ for any
$k\in\Z^N$ and $(t,x)\in\R^+\times\R^N$;
\item [(H4)]$u_0\in W^{2,\infty}(\R^N)$.
\end{itemize}

When $g\equiv C_N$, with $C_N$  a suitable constant depending on
the dimension $N$, then \eqref{levy} is the integral
representation of $-(-\Delta)^\frac{1}{2}$ for bounded real smooth
functions defined on $\R^N$ (see Theorem 1 in \cite{id}). We
recall that $(-\Delta)^\frac{1}{2}$ is the fractional operator
defined for instance on the Schwartz class  {\em S}$(\R^N)$ by
\beq\label{halflapfourier}\widehat{{(-\Delta )^\frac{1}{2}v}\
}(\xi)=|\xi|\ \widehat{v}(\xi),\eeq where $\widehat{w}$ is the
Fourier transform of $w$.

We prove that the limit $u^0$ of $u^\ep$ as $\ep\rightarrow0$
exists and is the unique solution of the homogenized problem
\begin{equation}\label{ueffett}
\begin{cases}
\p_{t} u=\overline{H}(\nabla_x u,\I[u(t,\cdot)])
&\text{in}\quad \R^+\times\R^N\\
u(0,x)=u_0(x)& \text{on}\quad \R^N,
\end{cases}
\end{equation}
for some continuous function $\overline{H}$ usually called {\em
effective Hamiltonian}.
The function $u^0$ will be interpreted later  as a macroscopic plastic strain
satisfying the macroscopic  plastic flow rule (\ref{ueffett}). 
Moreover $\I[u^0]$ will be the stress created by the macroscopic density of dislocations.

\subsection{Main results}
As usual in periodic homogenization, the limit equation is
determined by a {\em cell problem}. In our case, such a problem is
for any $p\in\R^N$ and $L\in \R$ the following:
\begin{equation}\label{v}
\begin{cases}
\lam+\p_{\tau} v=\I[v(\tau,\cdot)]+L-W'(v+\lam \tau+p\cdot y)+\s
(\tau,y)&\text{in}\quad \R^+\times\R^N\\ v(0,y)=0& \text{on}\quad
\R^N,
\end{cases}
\end{equation}
where $\lam=\lam(p,L)$ is the unique number for which there exists
a solution $v$ of \eqref{v} which is bounded on $\R^+\times\R^N$. In
order to solve \eqref{v}, we show for any $p\in\R^N$ and $L\in \R$
the existence of a unique solution of
\begin{equation}\label{w}
\begin{cases}
\p_{\tau} w=\I[w(\tau,\cdot)]+L-W'(w+p\cdot y)+\s
(\tau,y)&\text{in}\quad \R^+\times\R^N\\ w(0,y)=0& \text{on}\quad
\R^N,
\end{cases}
\end{equation}
and we look for some $\lam\in\R$ for which $w-\lam\tau$ is
bounded. Precisely we have:
\begin{thm}[Ergodicity]\label{ergodic}Assume (H1)-(H4). For $L\in\R$ and $p\in\R^N$, there
exists a unique viscosity solution $w\in C_b(\R^+\times\R^N)$ of
\eqref{w} and there exists a unique $\lam\in\R$ such that $w$
satisfies: $\frac{w(\tau,y)}{\tau}$ converges towards $\lam$ as
$\tau\rightarrow+\infty$, locally uniformly in $y$. The real
number  $\lam$ is denoted by $\overline{H}(p,L)$. The function
$\overline{H}(p,L)$ is continuous on $\R^N\times\R$ and
non-decreasing in $L$.
\end{thm}
 Unfortunately, we cannot directly use the bounded solution of \eqref{v}, usually called
{\em corrector}, in order to prove the convergence of the sequence
$u^\ep$ to the solution of \eqref{ueffett}. Nevertheless we have
the following result:
\begin{thm}[Convergence]\label{convergence}Assume (H1)-(H4). The solution
$u^\ep$ of \eqref{uep} converges towards the solution $u^0$ of
\eqref{ueffett} locally uniformly in $(t,x)$, where $\overline{H}$
is defined in Theorem \ref{ergodic}.
\end{thm}

Let us mention that in a companion paper \cite{mp}, we show that we can recover Orowan's law in dimension $N=1$ for $\sigma=0$, i.e.
$$\overline{H}(\delta p,\delta L) \simeq c_0 \delta^2 |p| L \quad \mbox{as}\quad \delta \to 0$$
i.e. the plastic strain velocity is asymptotically proportional to the product of dislocation density $|p|$ by the effective stress $L$.
\bigskip

\subsection{Brief review of the literature.}\label{strategyhom}

This non-local equation \eqref{uep} is related to the local
equation \beq\label{localintro}\begin{cases}
\p_tu^\ep =F\left(\frac{x}{\ep},\frac{u^\ep}{\ep},\nabla u^\ep\right)&\text{in}\quad \R^+\times\R^N\\
u^\epsilon(0,x)=u_0(x)& \text{on}\quad \R^N,\end{cases}\eeq that
 was studied in \cite{im} under the assumption that
$F(x,u,p)$ is periodic in $(x,u)$ and coercive in $p$. The
homogenization problem \eqref{localintro} when $F$ does not depend
on $u$, has been completely solved by Lions Papanicolaou and
Varadhan \cite{lpv}. After this seminal paper, homogenization of
Hamilton-Jacobi equations for coercive Hamiltonians has been
treated  for a wider class of periodic situations, c.f. Ishii
\cite{ishii}, for problems set on bounded domains, c.f. Alvarez
\cite{alv}, Horie and Ishii \cite{hi}, for equations with
different structures, c.f. Alvarez and Ishii \cite{ai}, for
deterministic control problems in $L^\infty$, c.f. Alvarez and
Barron \cite{ab}, for almost periodic Hamiltonians, c.f. Ishii
\cite{ishii2}, and for Hamiltonians with stochastic dependence,
c.f. Souganidis \cite{su}. More recently, inspired by \cite{im},
Barles \cite{b} gave an homogenization result for non-coercive
Hamiltonians and, as a by-product, obtained a simpler proof of the
results \cite{im} of Imbert and Monneau but under slightly more
restrictive assumptions on the Hamiltonians. We can also mention the work of Imbert,
Monneau and Rouy \cite{imr} where the authors studied homogenization of certain
integro-differential equations depending explicitly on
$u^\ep/\ep$. Notice that in the present paper, the operator $\I$ involves a singular kernel 
which creates some additional difficulties that were not present for instance in \cite{imr}.

Notice also that the model studied in \cite{imr} 
was introduced to approximate a level set model like in \cite{fim}.
The phase field model in \cite{imr} was therefore closer in the spirit 
to a model for discrete dislocation dynamics  at the mesoscopic scale.
On the contrary, the Peierls-Nabarro model (\ref{uep}) is a well-established physical model which
is really devoted to the description of dislocations at the microscopic scale.

\bigskip
\subsection{Organization of the paper} The paper is organized
as follows. In Section \ref{physicalsec}, we give more details
about the Peierls-Nabarro model yielding to the study of
\eqref{uep} and  the mechanical interpretation of the
homogenization results. In Section \ref{sect3} 
we present briefly the strategies of the main proofs.
In Section \ref{vissolsec}, we state
various com\-pa\-ri\-son principles, existence and regularity
results for solutions of non-local Hamilton-Jacobi equations. In
Section \ref{Proofconvsec}, we prove the convergence result
(Theorem \ref{convergence}) by assuming the existence of smooth
approximate sub and supercorrectors (Proposition
\ref{apprcorrectors}). In order to show their existence, in
Section \ref{lipcorrsec}, we first construct Lipschitz continuous
sub and supercorrectors (Proposition \ref{lipcorrect}). As a
byproduct,  we prove the ergodicity of the problem (Theorem
\ref{ergodic}) and some properties of the effective Hamiltonian
(Proposition \ref{Hprop}). Proposition \ref{apprcorrectors} is
then proved  in Section \ref{smoothcorsec}. 
The proofs of Lemma \ref{regularityvisc} 
and of Proposition \ref{pro::s19} are done in the Appendix (Section \ref{appendix}).

\subsection{Notations} We denote by $B_r(x)$ the ball of radius
$r$ centered at $x$. The cylinder $(t-\tau,t+\tau)\times B_r(x)$
is denoted by $Q_{\tau,r}(t,x)$.

$\lfloor x \rfloor$ and $\lceil x\rceil$ denote respectively the
floor and the ceil integer parts of a real number $x$.

It is convenient to introduce the singular measure defined on
$\R^N\setminus\{0\}$ by
$$\mu(dz)=\frac{1}{|z|^{N+1}}g\left(\frac{z}{|z|}\right)dz=\mu_0(z)dz,$$and
to denote
$$\I^{1,r}[U,x]=\int_{|z|\leq r}(U(x+z)-U(x)-\nabla U(x)\cdot
z)\mu(dz),$$
$$\I^{2,r}[U,x]=\int_{|z|>r}(U(x+z)-U(x))\mu(dz).$$
Sometimes when $r=1$ we will omit $r$ and we will write simply
$\I^1$ and $\I^2$.

For a function $u$ defined on $(0,T)\times\R^N$, $0<T\leq+\infty$,
for $0<\al<1$ we denote by $<u>_x^\al$  the seminorm defined by
$$<u>_x^\al:=\sup_{(t,x),\,(t,x')\in(0,T)\times\R^N\atop x\neq x'}\frac{|u(t,x)-u(t,x')|}{|x-x'|^\al}$$
and by $C_x^{\al}((0,T)\times\R^N)$ the space of continuous
functions defined on $(0,T)\times\R^N$ that are bounded and with
bounded seminorm
 $<u>_x^\al$.

Finally, we denote by $USC_b(\R^+\times\R^N)$ (resp.,
$LSC_b(\R^+\times\R^N)$) the set of upper (resp., lower)
semicontinuous functions on $\R^+\times\R^N$ which are bounded on
$(0,T)\times\R^N$ for any $T>0$ and we set
$C_b(\R^+\times\R^N):=USC_b(\R^+\times\R^N)\cap
LSC_b(\R^+\times\R^N)$.

\section{Physical modeling and mechanical interpretation of the homogenization results}\label{physicalsec}

\subsection{The Peierls-Nabarro model}
Dislocations are line defects in crystals. Their typical length is
of the order of  $10^{-6}m$ and their thickness of order of
$10^{-9}m$. When the material is submitted to shear stress, these
lines can move in the crystallographic planes and their dynamics
is one of the main explanation of the plastic behavior of metals.

The Peierls-Nabarro model is a phase field model for dislocation
dynamics incorporating atomic features into continuum framework.
In a phase field approach, the dislocations are represented by
transition of a continuous field.

We briefly review the model (see \cite{hl} for a detailed
presentation). As an example, consider an edge dislocation in a
crystal with simple cubic lattice.  In a Cartesian system of
coordinates $x_1x_2x_3$, we assume that the dislocation is located
in the slip plane $x_1x_2$ (where the dislocation can move) and
that the Burgers' vector  (i.e. a fixed vector associated to the
dislocation) is in the direction of the $x_1$ axis. We write this
Burgers' vector as $be_1$ for a real $b$. The disregistry of the
upper half crystal $\{x_3>0\}$ relative to the lower half
$\{x_3<0\}$ in the direction of the Burgers' vector is
$\phi(x_1,x_2)$, where $\phi$ is a phase parameter between $0$ and
$b$. Then the dislocation loop can be for instance localized by
the level set $\phi=b/2$. For a closed loop, we expect to have
$\phi\simeq b$ inside the loop and $\phi\simeq 0$ far outside the
loop.

In the Peierls-Nabarro model, the total energy is given by
\beq\label{energy}
\mathcal{E}=\mathcal{E}^{el}+\mathcal{E}^{mis}.\eeq

In \eqref{energy}, $\mathcal{E}^{mis}$ is the so called {\em misfit energy} due to
the nonlinear atomic interaction across the slip plane \beqs
\mathcal{E}^{mis}(\phi)=\int_{\R^2}W(\phi(x))\ dx\quad \mbox{with}\quad x=(x_1,x_2),\eeqs
where $W(\phi)$ is the interplanar potential. In the classical
Peierls-Nabarro model \cite{p,n}, $W(\phi)$ is
approximated by the
sinusoidal potential \beqs
W(\phi)=\frac{\mu
b^2}{4\pi^2d}\left(1-\cos\left(\frac{2\pi\phi}{b}\right)\right),\eeqs
where $d$ is the lattice spacing perpendicular to the slip plane.\\

The elastic energy $\mathcal{E}^{el}$ induced by the dislocation  is (for $X=(x,x_3)$ with $x=(x_1,x_2)$)
\beqs
\mathcal{E}^{el}(\phi,U)=\frac{1}{2}\int_{\R^3}e:\Lambda:e\
dX \quad \mbox{with}\quad e=e(U)-\phi(x)\delta_0(x_3) e^0 \quad \mbox{and}\quad \left\{\begin{array}{l}
e(U)=\frac12\left(\nabla U + (\nabla U)^T\right)\\
\\
e^0=\frac12\left(e_1\otimes e_3 + e_3\otimes e_1\right)
\end{array}\right.,\eeqs
where $U: \R^3\to \R^3$ is the displacement and
$\Lambda=\{\Lambda_{ijkl}\}$ are the elastic coefficients.\\
Given the field $\phi$, we minimize the energy $\mathcal{E}^{el}(\phi,U)$ with respect to the displacement $U$ and define
$$\mathcal{E}^{el}(\phi)=\inf_{U} \mathcal{E}^{el}(\phi,U)$$
Following the proof of Proposition 6.1 (iii) in \cite{ahlm}, we can see that (at least formally)
$$\mathcal{E}^{el}(\phi)=-\frac12 \int_{\R^2} (c_0\star\phi)\phi$$
where $c_0$ is a certain kernel.
In the case of isotropic elasticity, we have
$$\Lambda_{ijkl}=\lambda \delta_{ij}\delta_{kl} +\mu\left(\delta_{ik}\delta_{jl}+\delta_{il}\delta_{jk}\right)$$
where $\lambda,\mu$ are the Lam\'e coefficients. Then the kernel $c_0$ can be written
(see Proposition 6.2 in \cite{ahlm}, translated in our framework):
$$c_0(x)=\frac{\mu}{4\pi}\left(\partial_{22}\frac{1}{|x|}+ \gamma\partial_{11}\frac{1}{|x|}\right)\quad \mbox{with}\quad \gamma=\frac{1}{1-\nu}\quad \mbox{and}\quad \nu=\frac{\lambda}{2(\lambda+\mu)}$$
where $\nu\in (-1,1/2)$ is called the Poisson ratio.

The equilibrium configuration of straight dislocations is obtained by
minimizing the total energy with respect to $\phi$, under the constraint that far from the dislocation core, the function
 $\phi$ tends to $0$ in one half plane and to $b$ in the other half plane. In particular, the phase transition $\phi$ is then solution of the following equation
\beq\label{hallapleqn=1}\I[\phi]=W'(\phi)\quad\text{on }\R^2,\eeq
where formally $\I[\phi] = c_0\star \phi$, which is the anisotropic L\'evy operator defined in (\ref{levy}) for $N=2$
and $g(z_1,z_2)=\frac{\mu}{4\pi}\left((2\gamma-1)z_1^2+(2-\gamma)z_2^2\right)$.
Let us now recall the expression of the kernel after a Fourier transform (see paragraph 6.2.2.2 in \cite{ahlm})
$$\widehat{c_0}(\xi)=-\frac{\mu}{2|\xi|}\left(\xi_2^2+\gamma \xi_1^2\right)$$
Then for $\gamma=1$ and $\mu=2$, we see that $\I = -(-\Delta)^{\frac12}$.
In that special case, we recall that the solution $\phi$ of (\ref{hallapleqn=1}) satisfies $\phi(x)=\tilde{\phi}(x,0)$ where $\tilde{\phi}(X)$ is the solution of (see \cite{landkof,gonzalezmonneau})
$$\left\{\begin{array}{ll}
\Delta \tilde{\phi}=0 &\quad \mbox{in}\quad \left\{x_3>0\right\}\\
\\
\displaystyle \frac{\partial \tilde{\phi}}{\partial x_3} = W'( \tilde{\phi}) &\quad \mbox{on}\quad \left\{x_3=0\right\}
\end{array}\right.$$
Moreover, we have in particular an explicit solution for $b=1$, $d=2$ (with $W'(\tilde{\phi})=\frac1{2\pi}\sin (2\pi \tilde{\phi})$)
$$\displaystyle \tilde{\phi}(X)=\frac12 +\frac{1}{\pi}\arctan \left(\frac{x_1}{x_3+1}\right) $$

Then by rescaling, it is easy to check that we can recover the explicit solution found in Nabarro \cite{n}
\beqs
\left\{\begin{array}{ll}
\displaystyle \phi(x)=\frac{b}{2}+ \frac{b}{\pi}\arctan\left(\frac{2(1-\nu)x_1}{d}\right) & (\mbox{edge dislocation})\\
\\
\displaystyle \phi(x)=\frac{b}{2}+ \frac{b}{\pi}\arctan\left(\frac{2x_2}{d}\right) & (\mbox{screw dislocation})
\end{array}\right.\eeqs

In a more general model, one can consider a potential $W$ satisfying
\begin{itemize}
    \item [(i)] $W(v+b)=W(u)$ for all $v\in\R$;
    \item [(ii)]$W(b\Z)=0<W(a)$ for all $a\in\R\setminus b\Z$.
\end{itemize}
The periodicity of $W$ reflects the periodicity of the crystal,
while the minimum property is consistent with the fact that the
perfect crystal is assumed to minimize the energy.

\bigskip

In the face cubic structured (FCC) observed in many metals and
alloys, dislocations move at low temperature on the slip plane. In
the present paper we are interested in describing the effective
dynamics for a collection of dislocations curves with the same
Burgers' vector and all contained in a single slip plane $x_1x_2$,
and moving in a landscape with periodic obstacles (that can be for
instance precipitates in the material). These dislocations are
represented by a single phase parameter $u(t,x_1,x_2)$ defined on
the slip plane $x_1x_2$. The dynamic of dislocations  is then
described by the evolutive version  of the Peierls-Nabarro model
 (see for instance \cite{MBW} and \cite{Denoual}):
\begin{equation}\label{nabarroevolutintro}
\p_{t} u=\I[u(t,\cdot)]-W'\left(u\right)+\s_{13}^{\tiny\mbox{obst}}
\left(t,x\right)\quad\text{in}\quad \R^+\times\R^N\\
\end{equation}
for $x\in\R^N$ with the physical dimension $N=2$. In the model,
the component $\s_{13}^{\tiny\mbox{obst}}$ of the stress
(evaluated on the slip plane) has been introduced to take into
account the shear stress not created by the dislocations
themselves. This shear stress is created by the presence of the
periodic {\it obstacles} and the possible external applied stress
on the material.

We want to identify at {\em large scale} an evolution model for
the dynamics of a density of dislocations. We consider the
following rescaling
$$u^\ep(t,x)=\ep u\left(\frac{t}{\ep},\frac{x}{\ep}\right),$$
where $\epsilon$ is the ratio between the typical length scale for
dislocation (of the order of the micrometer) and the typical
macroscopic length scale in mechanics (milimeter or centimeter).
With such a rescaling, we see that the number of dislocations is
typically of the order of $1/\ep$ per unit of macroscopic scale.
Moreover, assuming suitable initial data
\beq\label{nabincondintro} u(0,x)=\frac{1}{\ep}u_0(\ep
x)\quad\text{on }\R^N,\eeq (where $u_0$ is a regular bounded
function), we see that the functions $u^\ep$ are solutions of
(\ref{uep}). This indicates that at the limit $\epsilon\to 0$, we
will recover a model for the dynamics of (renormalized) densities
of dislocations.

\begin{rem}{\em Fractional reaction-diffusion equations of the
form \beq\label{frcdiffintro}\p_t u=\I[u]+f(u)\quad\text{in
}\R^+\times\R^N\eeq where $N\geq 2$ and $f$ is a bistable
nonlinearity have been studied by Imbert and Souganidis \cite{is}.
In this paper the authors show that solutions of
\eqref{frcdiffintro}, after properly rescaling them, exhibit the limit evolution of an interface by (anisotropic)
mean curvature motion.

Other results have been obtained by
Gonz\'{a}lez and Monneau \cite{gonzalezmonneau} for a rescaling of the evolutive
Peierls-Nabarro model in dimension $N=1$.
In the one dimensional space, the limit moving interfaces are points particles interacting with forces as $1/x$.
The dynamics of these
particles corresponds to the classical discrete dislocation
dynamics, in the particular case of parallel straight edge
dislocation lines in the same slip plane with the same Burgers'
vector.
In \cite{fim}, considering another rescaling of the model of particles obtained in \cite{gonzalezmonneau},
the authors identify at large scale an evolution model for the dynamics of a density of
dislocations, that is analoguous to \eqref{ueffett}. In the present
paper, we directly deduce  the model \eqref{ueffett} at larger
scale from the Peierls-Nabarro model at smaller scale in any dimension $N\ge 1$. That way we remove the limitation to the dimension $N=1$ that appears in \cite{gonzalezmonneau}.

Finally, let us mention that in \cite{gm} and \cite{gm2} Garroni
and Muller study a variational model for dislocations that is the
variational formulation of the stationary Peierls-Nabarro
equation, where they derive a line tension model.}\end{rem}

\subsection{Mechanical interpretation of the homogenization}
Let us briefly explain the meaning of the homogenization result.
In the macroscopic model, the function $u^0(t,x)$
can be interpreted as the plastic strain (localized in the slip plane $\left\{x_3=0\right\}$).
Then the three-dimensional displacement $U(t,X)$  is obtained as a minimizer of the elastic energy
$$\displaystyle U(t,\cdot)={\mbox{arg}\min_{\tilde{U}}}\ \mathcal{E}^{el}(u^0(t,\cdot),\tilde{U})$$
and the stress is
$$\sigma=\Lambda: e \quad \mbox{with}\quad e=e(U)-u^0(t,x)\delta_0(x_3)e^0$$
Then the resolved shear stress is
$$\I[u^0]= \sigma^{\tiny\mbox{obst}}_{13}$$
The homogenized equation \eqref{ueffett}, i.e.
$$\p_{t} u^0=\overline{H}(\nabla_x u^0,\I[u^0(t,\cdot)])$$
which is the evolution equation for $u^0$, can be interpreted as the plastic flow rule in a model for macroscopic crystal plasticity.
This is the law giving the plastic strain velocity $\p_t u^0$ as a function of the resolved shear stress
$\sigma^{\tiny\mbox{obst}}_{13}$ and the dislocation density $\nabla u^0$.

The typical example of such a plastic flow rule is the Orowan's law:
$$\overline{H}(p,L) \simeq  |p|L$$
This is also the law that we recover in dimension $N=1$ in a forthcoming paper \cite{mp}
 in the case where there are no obstacles (i.e.
$\sigma_{13}^{\tiny\mbox{obst}}\equiv 0$) and for small stress $L$ and small density $|p|$.
When $\sigma_{13}^{\tiny\mbox{obst}}\not\equiv 0$ with zero mean value  (i.e. $< \sigma^{\tiny\mbox{obst}}_{13}> =0$), 
we expect a threshold phenomenon as in \cite{imr} (see also
Norton's law with threshold in \cite{FPZ}), i.e.
$$\overline{H}(p,L)=0\quad \mbox{if}\quad |L|\quad \mbox{is small enough}.$$
This means more generally that our homogenization procedure
describes correctly the mechanical behaviour of the stress at
large scales, but keeps the memory of the microstructure in the
plastic law with possible threshold effects.



\section{Strategies of the main proofs}\label{sect3}

\subsection{Strategy for the proof of convergence}{$\mbox{ }$}
\subsubsection{The general approach}{$\mbox{ }$}
It has been already noticed that for problems periodic in $u^\ep/\ep$, we have to introduce twisted correctors
(see for instance \cite{im}). It is also known that if we can claim that the limit function satisfies
\begin{equation}\label{eq::s16}
\partial_t u^0 \not= 0 \quad \mbox{or}\quad \nabla_x u^0\not= 0
\end{equation}
then we do not have to introduce an additional dimension to perform the proof of convergence. The idea (see \cite{im}) is that we can twist the corrector either dividing by $p_i:=\partial_{x^i} u^0$ for some index $i$, or by $\lambda:=\partial_t u^0$ like considering the ansatz:
$$u^\varepsilon(t,x) \simeq u^0(t,x) + \varepsilon v\left(\frac{u^0(t,x)-p\cdot x}{\varepsilon \lambda}, \frac{x}{\varepsilon}\right).$$
On the contrary, we do not know how to deal with the case where both quantities in (\ref{eq::s16})
vanish, except adding a dimension and considering twisted correctors in higher dimension.
Here we have to face a
similar difficulty in the much more involved framework of
non-local equations. Notice also  that it does not seem possible to apply the
approach of Barles \cite{b}. Therefore following the idea in
\cite{im}, we consider the solution $U^\ep$ of
 \begin{equation}\label{Uepintro}
\begin{cases}
\p_{t}
U^\epsilon=\I[U^\epsilon(t,\cdot,x_{N+1})]-W'\left(\frac{U^\epsilon}{\epsilon}\right)+\s
\left(\frac{t}{\epsilon},\frac{x}{\epsilon}\right)&\text{in}\quad \R^+\times\R^{N+1}\\
U^\epsilon(0,x,x_{N+1})=u_0(x)+p_{N+1}x_{N+1}& \text{on}\quad
\R^{N+1},
\end{cases}
\end{equation}where $p_{N+1}\neq0$. We then consider the following ansatz:
\beqs U^\ep(t,x,x_{N+1})\simeq U^0(t,x,x_{N+1})+\ep
V\left(\frac{t}{\ep},\frac{x}{\ep},\frac{U^0(t,x,x_{N+1})-\lam
t-p\cdot x}{\ep p_{N+1} }\right)\eeqs
 where $U^0 (t,x,x_{N+1})=u^0(t,x)+ p_{N+1}x_{N+1}$. This ansatz turns out to be the good one, and plugging this expression of $U^\ep$
 into \eqref{Uepintro}, we find formally with
 $\tau=\frac{t}{\ep},\,y=\frac{x}{\ep}$, $y_{N+1}=\frac{U^0(t,x,x_{N+1})-\lam t-p\cdot
x}{p_{N+1}\ep }$:
 \beq\label{Vintro}\lam+\p_{\tau} V=L+\I[V(\tau,\cdot,y_{N+1})]-W'(V+p\cdot y+p_{N+1}y_{N+1}+\lam\tau)+\s(\tau,y),\eeq
     where $$\lam=\p_t U^0(t,x,x_{N+1})=\p_t u^0(t,x),\quad p=\nabla_x U^0(t,x,x_{N+1})=\nabla_x u^0(t,x)$$ and
     $$L=\I[U^0(t,\cdot,x_{N+1})]=\I[u^0(t,\cdot)].$$  
     Then, we expect $u^0$ to be solution of \eqref{ueffett} with $\bar{H}(p,L)=\lam(p,L)$. 
     This heuristic computation, that  permits first of all to identify the cell problem
     in the higher dimensional space, can be made rigorous through the perturbed test function method by Evans \cite{e1}.\\

\subsubsection{Additional difficulty}{$\mbox{ }$}\label{sect312}

     Let us enter a bit more in the details of the proof.  Fix $P_0=(t_0,x_0,x_{N+1}^0)\in\R^+\times\R^{N+1}$ and define
     \beq\label{tildeUep} \tilde{U}^\ep(t,x,x_{N+1})= U^0(t,x,x_{N+1})+\ep
V\left(\frac{t}{\ep},\frac{x}{\ep},\frac{U^0(t,x,x_{N+1})-\lam
t-p\cdot x}{\ep p_{N+1} }\right),\eeq where $V$ is solution of \eqref{Vintro} with $\lam=\p_t U^0(P_0),$ $ p=\nabla_x U^0(P_0)$ and
     $L=\I[U^0(t_0,\cdot,x_{N+1}^0),x_0]$. Let us call  $F(t,x,x_{N+1})=\frac{U^0(t,x,x_{N+1})-\lam
t-p\cdot x}{ p_{N+1} }$. Here we assume for simplicity that $U^0$ and $V$ are smooth. 
The proof of convergence consists in showing that  $\tilde{U}^\ep$ is a solution of \eqref{Uepintro} in a cylinder  $(t_0-r,t_0+r)\times B_r(x_0,x_{N+1}^0)$ for $r>0$ small enough, up to an error that goes to 0 as $r\rightarrow 0^+$.  This will allow us to compare $U^\epsilon$ with $\tilde{U}^\ep$ and, thanks to the boundedness of $V$, to conclude that $U^\epsilon$ converges to $U^0$ as $\ep\rightarrow 0$.  

When we plug $\tilde{U}^\ep$ into \eqref{Uepintro},  we find the equation 
\beqs \lam + \p_{\tau} V=L+\I[V(\tau,\cdot,y_{N+1})]-W'(V+p\cdot y+p_{N+1}y_{N+1}+\lam\tau)+\s(\tau,y)+o_r(1)+\theta_r,\eeqs

with   $\tau=\frac{t}{\ep},\,y=\frac{x}{\ep}$, $y_{N+1}=\frac{F(t,x,x_{N+1})}{\ep }$, where 
\beqs   \begin{split} \theta_r &= (\p_tU^0(P_0)-\p_t U^0(t,x,x_{N+1}))\p_{y_{N+1}}V(\tau,y,y_{N+1})
\\&+\I\left[V\left(\tau,\cdot, \frac{F(\ep\tau,\ep\cdot, \ep y_{N+1})}{\ep}\right)\right]-\I[V(\tau,\cdot,y_{N+1})].\end{split} \eeqs

Then,       $\tilde{U}^\ep$ will be a solution of  \eqref{Uepintro} up to a small error if $\theta_r=o_r(1)$ as $r\rightarrow0+$.  This last property holds true if the corrector $V$  satisfies: 
$|V|$, $|\p_{y_{N+1}}V|\leq C$ in $\R^+\times\R^{N+1}$ for some $C>0$, and 
\begin{equation}\label{eq::s17}
\p_{y_{N+1}}V(\tau,\cdot,\cdot) \quad \mbox{is H\"{o}lder continuous, uniformly in time.}
\end{equation}
 In the case of the local first order equation \eqref{localintro} considered  in \cite{im},
or non local equations considered in \cite{imr},  
approximate correctors were only required to be Lipschitz continuous in the additional variable. 
Here the additional regularity (\ref{eq::s17}) is required because we deal with an operator $\I$ whose kernel is singular.

Since in \eqref{Vintro}, the quantity $\I[V(\tau,\cdot,y_{N+1})]$ is computed only in the $y$ variable, we cannot expect this kind of regularity
 for the correctors.
       Nevertheless, we are able to construct regular
 approximated sub and supercorrectors, i.e., sub and supersolutions of approximate $N+1$-dimensional cell problems, and this is enough to conclude. 
  Finally, this construction works for any $p_{N+1}\neq0$ and to simplify the presentation we take $p_{N+1}=1$.

 \bigskip
 
\subsection{Strategy for the construction of smooth approximate correctors}\label{sect3.2}

As explained in the previous subsection, in the proof of convergence we will need smooth  approximate sub and 
and super-correctors on $\R^+\times\R^{N+1}$, i.e.,   for $P=(p,1)\in\R^{N+1}$ and $L\in\R$, sub and supersolutions of
\begin{equation}\label{V}\left\{
  \begin{array}{ll}
    \lam+\p_{\tau} V=L+\I[V(\tau,\cdot,y_{N+1})]-W'(V+P\cdot Y+\lam\tau)+\s(\tau,y) & \hbox{in } \R^+\times\R^{N+1}\\
    V(0,Y)=0 & \hbox{on }\R^{N+1}.
  \end{array}
\right.\end{equation} Here and in what follows, we denote
$Y=(y,y_{N+1})$. 
More precisely, we will prove the following proposition.
\begin{prop}[Smooth approximate correctors]\label{apprcorrectors}
Let $\lam$ be the constant defined by Theorem \ref{ergodic}. For
any fixed $p\in\R^N$, $P=(p,1)$, $L\in\R$ and $\eta>0$ small
enough, there exist real numbers $\lam^+_\eta(p,L)$,
$\lam^-_\eta(p,L)$, a constant $C>0$ (independent of $\eta,\,p$
and $L$) and bounded super and subcorrectors $V^+_{\eta},
V^-_{\eta}$, i.e. respectively a super and a subsolution of
\begin{equation}\label{apprcorrequ}\left\{
  \begin{array}{ll}
    \lam^{\pm}_\eta+\p_{\tau}
    V^{\pm}_\eta=L+\I[V^{\pm}_\eta(\tau,\cdot,y_{N+1})]\\
    \qquad\qquad\qquad-W'(V^{\pm}_\eta+P\cdot Y+\lam^{\pm}_\eta\tau)+\s(\tau,y){\mp} o_\eta(1) & \hbox{in } \R^+\times\R^{N+1}\\
    V^{\pm}_\eta(0,Y)=0 & \hbox{on }\R^{N+1},
  \end{array}
\right.\end{equation}  where $0\leq o_\eta(1)\rightarrow0$ as
$\eta\rightarrow0^+$, such that \beq\label{appcorr1}
\lim_{\eta\rightarrow0^+}\lam^+_\eta(p,L)=\lim_{\eta\rightarrow0^+}\lam^-_\eta(p,L)=\lam(p,L),\eeq
locally uniformly in $(p,L)$, $\lam^{\pm}_\eta$ satisfy (i) and
(ii) of Proposition \ref{Hprop} and for any
$(\tau,Y)\in\R^+\times\R^{N+1}$
\beq\label{appcorr2}|V^{\pm}_{\eta}(\tau,Y)|\leq C.\eeq
 Moreover $V^{\pm}_{\eta}$ are of class $C^{2}$ w.r.t. $y_{N+1}$, and for any $0<\al<1$ 
\beq\label{appcorr3}
-1\leq\p_{y_{N+1}}V^{\pm}_\eta \leq
\frac{\|W''\|_\infty}{\eta},\eeq 
\begin{equation}\label{contrderivappcorr}
\|\p^2_{y_{N+1}y_{N+1}}V^{\pm}_\eta\|_\infty\leq C_{\eta},\quad <\p_{y_{N+1}}V^{\pm}_\eta>_y^\al,\,\le C_{\eta,\alpha}.
\end{equation}
\end{prop}

Here in order to build Lipschitz sub/super correctors, 
it does not seem easy to apply a kind of truncation of the Hamiltonian
like in \cite{im} or \cite{imr}. Therefore we use a 
different method to build such approximate correctors (similar to the one in \cite{fim2}).

The proof of  Proposition \ref{apprcorrectors} is mainly performed in two steps:

\noindent {\bf Step 1: Constructions of Lipschitz correctors.}\\
Using the modified Cauchy problem
$$\left\{\begin{array}{llll}
\partial_\tau U &= &L + \I[U(\tau,\cdot,y_{N+1})] - W'(U+P\cdot Y) +\sigma(\tau,y)&\\
\\ & &+\eta \left\{a_0 +\displaystyle \inf_{Y'} U(\tau,Y')-U(\tau,Y)\right\}|\partial_{y_{N+1}}U+1| & \quad \mbox{in}\quad \R^+\times \R^{N+1}\\
\\
U(0,Y)&=& 0 & \quad \mbox{on}\quad \R^{N+1}
\end{array}\right.$$
we construct Lipschitz correctors.
The Lipschitz bound comes formally from the equation satisfied by $w=\partial_{y_{N+1}}U$:
$$\left\{\begin{array}{llll}
\partial_\tau w &= & \I[w(\tau,\cdot,y_{N+1})] - W''(U+P\cdot Y)(w+1)-\eta w(\tau,Y)|w+1|&\\
\\ & &+  \eta \left\{a_0 +\displaystyle \inf_{Y'} U(\tau,Y')-U(\tau,Y)\right\}\ \mbox{sign}(\partial_{y_{N+1}}U+1)
 \partial_{y_{N+1}}w
& \quad \mbox{in}\quad \R^+\times \R^{N+1}\\
\\
w(0,Y)&=& 0 & \quad \mbox{on}\quad \R^{N+1}
\end{array}\right.$$
and the comparison principle implies that
\begin{equation}\label{eq::s18}
-1\le w\le \frac{|W''|_\infty}{\eta}
\end{equation}
On the other hand we are able to show (as in \cite{imr}) that $\inf_{Y'} U(\tau,Y')-U(\tau,Y)$ remains bounded independently on $\eta$.
Then an appropriate choice of $a_0$ large enough (resp.  negative enough) provides us bounded supercorrectors $W^+_\eta$ 
(resp.  subcorrectors $W^-_\eta$). We also show using Proposition \ref{regularityvisc} and the bound (\ref{eq::s18})
that we have the following H\"{o}lder estimate:
$$<W^\pm_\eta >^\alpha_y \quad \le \quad  C_\alpha$$


\noindent {\bf Step 2: Constructions of smooth correctors.}\\
We make a convolution with respect to $y_{N+1}$ of the Lipschitz correctors built in Step 1,  with a sequence  $(\rho_\delta)_\delta$ of mollifiers:
\beqs V^{\pm}_{\eta,\delta}(t,y,y_{N+1}):=W^{\pm}_\eta(t,y,\cdot)\star\rho_\delta(\cdot).\eeqs
Those functions are finally the smooth approximate sub/super correctors of Proposition \ref{apprcorrectors}
with some small error term $o_\eta(1)$ on the right hand side of the equation, for a suitable choice $\delta=\delta(\eta)$.

\section{Results about viscosity solutions for non-local equations}\label{vissolsec}
The classical notion of viscosity solution can be adapted for
Hamilton-Jacobi equations involving non-local operators, see for
instance \cite{s}. In this section we state comparison principles,
existence and regularity results for viscosity solutions of
\eqref{uep} and \eqref{ueffett}, that will be used later in the
proofs.
\subsection{Definition of viscosity solution}
We first recall the definition of viscosity solution for a general
first order non-local equation with associated initial
condition:
\begin{equation}\label{generalpb}
\begin{cases}
u_t=F(t,x,u,Du,\I[u])&\text{in}\quad \R^+\times\R^N\\
u(0,x)=u_0(x)& \text{on}\quad \R^N,
\end{cases}
\end{equation}where $F(t,x,u,p,L)$ is continuous and
non-decreasing in $L$.
\begin{de}[r-viscosity solution]\label{defviscosity}A function $u\in USC_b(\R^+\times\R^N)$ (resp., $u\in LSC_b(\R^+\times\R^N)$) is a
$r$-viscosity subsolution (resp., supersolution) of
\eqref{generalpb} if $u(0,x)\leq (u_0)^*(x)$ (resp., $u(0,x)\geq
(u_0)_*(x)$) and for any $(t_0,x_0)\in\R^+\times\R^N$, any
$\tau\in(0,t_0)$ and any test function $\phi\in
C^2(\R^+\times\R^N)$ such that $u-\phi$ attains a local maximum
(resp., minimum) at the point $(t_0,x_0)$ on
$Q_{(\tau,r)}(t_0,x_0)$, then we have
\beqs\begin{split}&\p_t\phi(t_0,x_0)-F(t_0,x_0,u(t_0,x_0),\nabla_x\phi(t_0,x_0),\I^{1,r}[\phi(t_0,\cdot),x_0]+\I^{2,r}[u(t_0,\cdot),x_0])\leq
0\\&\text{(resp., }\geq 0).\end{split}\eeqs  A function $u\in
C_b(\R^+\times\R^N)$ is a $r$-viscosity solution of
\eqref{generalpb} if it is a $r$-viscosity sub and supersolution
of \eqref{generalpb}.
\end{de}

It is classical that the maximum in the above definition can be
supposed to be global and this will be used later.  We have also
the following property, see e.g. \cite{s}:
\begin{prop}[Equivalence of the definitions]Assume $F(t,x,u,p,L)$ conti\-nuous and non-decreasing in $L$. Let $r>0$ and $r'>0$.
A function $u\in USC_b(\R^+\times\R^N)$ (resp., $u\in
LSC_b(\R^+\times\R^N)$) is a $r$-viscosity subsolution (resp.,
supersolution) of \eqref{generalpb} if and only if it is a
$r'$-viscosity subsolution (resp., supersolution) of
\eqref{generalpb}.
\end{prop}
Because of this proposition, if we do not need to emphasize $r$,
we will omit it when calling viscosity sub and supersolutions.
\subsection{Comparison principle and existence results}
In this subsection, we successively give com\-pa\-ri\-son
principles and existence results for \eqref{uep} and
\eqref{ueffett}. The following comparison theorem is shown in
\cite{jk} for more general parabolic integro-PDEs.
\begin{prop}[Comparison Principle for \eqref{uep}]\label{comparisonuep} Consider
 $u\in USC_b(\R^+\times\R^N)$ subsolution
and $v\in LSC_b(\R^+\times\R^N)$ supersolution of \eqref{uep},
then $u\leq v$ on $\R^+\times\R^N$.
\end{prop}
Following \cite{jk} it can  also be proved the comparison
principle for \eqref{uep} in bounded domains. Since we deal with a
non-local equation, we need to compare the sub and the
supersolution everywhere outside the domain.
\begin{prop}[Comparison Principle on bounded domains for
\eqref{uep}]\label{comparisonbounded} Let $\Om$ be a bounded
domain of $\R^+\times\R^N$ and let $u\in USC_b(\R^+\times\R^N)$
and $v\in LSC_b(\R^+\times\R^N)$ be respectively a sub and a
supersolution of $$\p_{t}
u^\epsilon=\I[u^\epsilon(t,\cdot)]-W'\left(\frac{u^\epsilon}{\epsilon}\right)+\s
\left(\frac{t}{\epsilon},\frac{x}{\epsilon}\right)$$ in $\Om$. If
$u\leq v$ outside $\Om$, then $u\leq v$ in $\Om$.
\end{prop}

\begin{prop}[Existence for \eqref{uep}]\label{existuep}For $\ep>0$ there exists
$u^{\ep}\in C_b(\R^+\times\R^N)$ (unique) viscosity solution of
\eqref{uep}. Moreover, there exists a constant $C>0$ independent
of $\ep$ such that
\begin{equation}\label{stimathmexist}|u^\ep(t,x)-u_0(x)|\leq
Ct.\end{equation}
\end{prop}
\dim Adapting the argument of \cite{imbert}, we can construct a
solution by Perron's method if we construct sub and supersolutions
of \eqref{uep}. Since $u_0\in W^{2,\infty}$, the two functions
$u^{\pm}(t,x):=u_0(x){\pm} Ct$ are respectively a super and a
subsolution of \eqref{uep} for any $\ep>0$, if
$$C\geq D_N\|u_0\|_{2,\infty}+\|W'\|_\infty+\|\sigma\|_\infty,$$
with $D_N$ depending on the dimension $N$. By comparison we also
get the estimate \eqref{stimathmexist}.\finedim We next recall the
comparison and the existence results for \eqref{ueffett}.
\begin{prop}[\cite{imr}, Proposition 3]\label{existHeff}
Let $\overline{H}: \R^N\times\R\rightarrow\R$ be continuous with
$\overline{H}(p,\cdot)$ non-decreasing on $\R$ for any $p\in\R^N$.
If $u\in USC_b(\R^+\times\R^N)$ and  $v\in LSC_b(\R^+\times\R^N)$
are respectively a sub and a supersolution of \eqref{ueffett},
then $u\leq v$ on $\R^+\times\R^N$. Moreover there exists a
(unique) viscosity solution of \eqref{ueffett}.
\end{prop}

In the next sections, we will embed the problem in the higher
dimensional space $\R^+\times\R^{N+1}$ by adding a new variable
$x_{N+1}$ in the equations. We will need the following proposition
showing that sub and supersolutions of the higher dimensional
problem are also sub and supersolutions  of the lower dimensional
one. This in particular implies that the comparison principle
between sub and supersolutions remains true increasing the
dimension.

\begin{prop}\label{fromN+1toN} Assume $F(t,x,x_{N+1},U,p,L)$ continuous and
 non-decreasing in $L$. Suppose that $U\in LSC_b(\R^+\times\R^{N+1})$ (resp.,  $U\in USC_b(\R^+\times\R^{N+1})$) is a viscosity
supersolution (resp., subsolution) of
\begin{equation}\label{eq::s10}
U_t=F(t,x,x_{N+1},U,D_xU,\I[U(t,\cdot,x_{N+1})])\quad\text{in}\quad
\R^+\times\R^{N+1},
\end{equation} 
then, for any $x_{N+1}\in\R$, $U$ is a
viscosity supersolution (resp., subsolution) of
$$U_t=F(t,x,x_{N+1},U,D_xU,\I[U(t,\cdot,x_{N+1})])\quad\text{in}\quad
\R^+\times\R^{N}.$$
\end{prop}
\dim 
Notice that in (\ref{eq::s10}), there is no derivative with respect to $x_{N+1}$ and no integral with respect to $dx_{N+1}$. Therefore $x_{N+1}$ only appears as a parameter that can (at least formally) be frozen.\\
We now do the (rigorous) proof for supersolutions. 
Fix $x^0_{N+1}\in\R$.
Let us consider a point  $(t_0,x_0)\in\R^+\times\R^N$ and a smooth
function $\varphi:\R^+\times\R^N\rightarrow\R$ such that
$$U(t,x,x^0_{N+1})-\varphi(t,x)\geq U(t_0,x_0,x^0_{N+1})-\varphi(t_0,x_0)=0\quad\text{for } (t,x)\in
Q_{\tau,r}(t_0,x_0),$$ with $r=1$. We have to show that
\beqs\begin{split}\p_t\varphi(t_0,x_0)&\geq
F(t_0,x_0,x^0_{N+1},U(t_0,x_0,x^0_{N+1}),D_x\varphi(t_0,x_0),\I^{1}[\varphi(t_0,\cdot),x_0]\\&+\I^{2}[U(t_0,\cdot,x^0_{N+1}),x_0]).\end{split}\eeqs
Without loss of generality, we can assume that the minimum is
strict. For $\ep>0$ let
$\varphi_\ep:\R^+\times\R^{N+1}\rightarrow\R$ be defined by
$$\varphi_\ep(t,x,x_{N+1})=\varphi(t,x)-\frac{1}{\ep}|x_{N+1}-x^0_{N+1}|^2.$$
Let $(t_\ep,x_\ep,x_{N+1}^\ep)$ be a minimum point of
$U-\varphi_\ep$ in $Q_{\tau,r}(t_0,x_0,x_{N+1}^0)$. Standard
arguments show that $(t_\ep,x_\ep,x_{N+1}^\ep)\rightarrow
(t_0,x_0,x_{N+1}^0)$ as $\ep\rightarrow0$ and that
$\lim_{\ep\rightarrow0}U(t_\ep,x_\ep,x_{N+1}^\ep)=U(t_0,x_0,x_{N+1}^0)$.
In particular, $(t_\ep,x_\ep,x_{N+1}^\ep)$ is internal to
$Q_{\tau,r}(t_0,x_0,x_{N+1}^0)$ for $\ep$ small enough, then we
get \beq\label{testrn+1}\p_t\varphi(t_\ep,x_\ep)\geq
F(t_\ep,x_\ep,U(t_\ep,x_\ep,x^\ep_{N+1}),D_x\varphi(t_\ep,x_\ep),\I^{1}[\varphi(t_\ep,\cdot),x_\ep]+\I^{2}[U(t_\ep,\cdot,x^\ep_{N+1}),x_\ep]).\eeq
By the Dominate Convergence Theorem
$\lim_{\ep\rightarrow0}\I^{1}[\varphi(t_\ep,\cdot),x_\ep]=\I^{1}[\varphi(t_0,\cdot),x_0]$;
by the Fatou's Lemma and  the convergence of
$U(t_\ep,x_\ep,x_{N+1}^\ep)$ to $U(t_0,x_0,x_{N+1}^0)$, we deduce
that
$$\I^{2}[U(t_0,\cdot,x^0_{N+1}),x_0]\leq\liminf_{\ep\rightarrow0}\I^{2}[U(t_\ep,\cdot,x^\ep_{N+1}),x_\ep].$$
Then, passing to the limit in \eqref{testrn+1} and using the
continuity and monotonicity of $F$, we get the desired inequality.
\finedim

\subsection{H\"{o}lder regularity}
In this subsection we state  a regularity result for sub
and supersolutions of semilinear non-local equations.
The proof is postponed in the appendix.
\begin{prop}[H\"{o}lder regularity]\label{regularityvisc}
Assume (H1) and let $g_1,\,g_2\in\R$. Suppose that $u\in
C(\R^+\times\R^{N})$ and bounded on $\R^+\times\R^{N}$ is a
viscosity subsolution of
\begin{equation*}
\begin{cases}
\p_{t} u=\I[u(t,\cdot)]+g_1&\text{in}\quad \R^+\times\R^N\\
u(0,x)=0& \text{on}\quad \R^N,
\end{cases}
\end{equation*} and a viscosity supersolution of
\begin{equation*}
\begin{cases}
\p_{t} u=\I[u(t,\cdot)]+g_2&\text{in}\quad \R^+\times\R^N\\
u(0,x)=0& \text{on}\quad \R^N.
\end{cases}
\end{equation*}
 Then, for any $0<\al<1$, $u\in
C_x^{\al}(\R^+\times\R^{N})$ with $<u>_x^\al\leq C$, where $C$
depends on $\|u\|_\infty,\,g_1$ and $g_2$.
\end{prop}

Notice that this regularity result will be used to establish a bound on the H\"{o}lder regularity in $y$ 
of $\partial_{y_{N+1}}V_\eta^{\pm}$ for smooth approximate correctors $V_\eta^{\pm}$ 
that will be used in Step 1.2 of the proof of Lemma \ref{converglem} 
used in the proof of the convergence result (Theorem \ref{convergence}).

\section{The proof of convergence}\label{Proofconvsec}

This section is dedicated to the proof of Theorem
\ref{convergence}. As explained in Subsection \ref{strategyhom},  we imbed our
problem in a higher dimensional one. We consider
$U^\ep$ solution of
\begin{equation}\label{Uep}
\begin{cases}
\p_{t}
U^\epsilon=\I[U^\epsilon(t,\cdot,x_{N+1})]-W'\left(\frac{U^\epsilon}{\epsilon}\right)+\s
\left(\frac{t}{\epsilon},\frac{x}{\epsilon}\right)&\text{in}\quad \R^+\times\R^{N+1}\\
U^\epsilon(0,x,x_{N+1})=u_0(x)+x_{N+1}& \text{on}\quad \R^{N+1}.
\end{cases}
\end{equation}
By Proposition \ref{fromN+1toN} and Proposition
\ref{comparisonuep}, the comparison principle holds true for
\eqref{Uep}. Then, as in the proof of Proposition  \ref{existuep},
by Perron's method we have:
\begin{prop}[Existence for \eqref{Uep}]\label{existUep}For $\ep>0$ there exists
$U^{\ep}\in C_b(\R^+\times\R^{N+1})$ (unique) viscosity solution
of \eqref{Uep}. Moreover, there exists a constant $C>0$
independent of $\ep$ such that
\begin{equation}\label{stimathmexistUep}|U^\ep(t,x,x_{N+1})-u_0(x)-x_{N+1}|\leq
Ct.\end{equation}
\end{prop}
 Let us exhibit the link between the
problem in $\R^N$ and the problem in $\R^{N+1}$.

\begin{lem}[Link between the problems on $\R^N$ and on $\R^{N+1}$]\label{linkuepUep} If $u^\ep$ and $U^\ep$ denote respectively the solution
of \eqref{uep} and \eqref{Uep}, then we have
$$\left|U^\ep(t,x,x_{N+1})-u^\ep(t,x)-\ep\left\lfloor\frac{x_{N+1}}{\ep}\right\rfloor\right|\leq
\ep,$$
\beq\label{linklemmdis}U^\ep\left(t,x,x_{N+1}+\ep\left\lfloor\frac{a}{\ep}\right\rfloor\right)
=U^\ep(t,x,x_{N+1})+\ep\left\lfloor\frac{a}{\ep}\right\rfloor\quad\text{for
any }a\in\R.\eeq
\end{lem}
This lemma is a consequence of the comparison principle for
\eqref{Uep}, the invariance by $\ep$-translations w.r.t. $x_{N+1}$
and the monotonicity of $U^\ep$ w.r.t. $x_{N+1}$.

 Let us now consider the problem
\begin{equation}\label{Ueffett}
\begin{cases}
\p_{t} U=\overline{H}(\nabla_x U,\I[U(t,\cdot,x_{N+1})])
&\text{in}\quad \R^+\times\R^{N+1}\\
U(0,x,x_{N+1})=u_0(x)+x_{N+1}& \text{on}\quad \R^{N+1}.
\end{cases}
\end{equation}
The link between problems \eqref{ueffett} and \eqref{Ueffett} is
given by the following lemma (analogue to Lemma \ref{linkuepUep}).
\begin{lem}\label{linkuU} Let $u^0$ and $U^0$ be respectively the solutions of \eqref{ueffett} and \eqref{Ueffett}. Then, we have
$$U^0(t,x,x_{N+1})=u^0(t,x)+x_{N+1}.$$
\end{lem}
Lemma \ref{linkuU} is a consequence of the comparison principle for
\eqref{Ueffett}  and the invariance by translations w.r.t.
$x_{N+1}$.

We need to make more precise the dependence of the real number
$\lam$ given by Theorem \ref{ergodic} on its variables. The
following properties will be shown in the next section.
\begin{prop}[Properties of the effective Hamiltonian]\label{Hprop} Let $p\in\R^N$ and $L\in\R$. Let $\overline{H}(p,L)$ be the constant
defined by Theorem \ref{ergodic}, then
$\overline{H}:\R^N\times\R\rightarrow\R$ is a continuous function
with the following properties:
\begin{itemize}
\item[(i)]$\overline{H}(p,L)\rightarrow {\pm}\infty$ as $L \rightarrow
{\pm}\infty$ for any $p\in\R^N$;
\item[(ii)] $\overline{H}(p,\cdot)$ is non-decreasing on $\R$ for any $p\in\R^N$;
\item[(iii)]If $\s(\tau,y)=\s(\tau,-y)$ then
$$\overline{H}(p,L)=\overline{H}(-p,L);$$
\item[(iv)]If $W'(-s)=-W'(s)$ and $\s(\tau,-y)=-\s(\tau,y)$ then
$$\overline{H}(p,-L)=-\overline{H}(p,L).$$
\end{itemize}
\end{prop}

\subsection{Proof of Theorem \ref{convergence}}{$\mbox{ }$}\\
\noindent {\bf Step 1: The classical approach}\\
By \eqref{stimathmexistUep}, we know that the
family of functions $\{U^\ep\}_{\ep>0}$ is locally bounded, then
$U^+:=\limsup^*_{\ep\rightarrow0}U^\epsilon$ is everywhere finite.
Classically we prove that
 $U^+$ is a subsolution of (\ref{Ueffett}).

Similarly, we can prove that
$U^-={\liminf_*}_{\ep\rightarrow0}U^\epsilon$ is a supersolution
of \eqref{Ueffett}. Moreover
$U^+(0,x,x_{N+1})=U^-(0,x,x_{N+1})=u_0(x)+x_{N+1}$. The comparison
principle for \eqref{Ueffett}, which is an immediate consequence
of Propositions \ref{existHeff} and \ref{fromN+1toN},
 then  implies that $U^+\leq U^-$. Since the reverse inequality $U^-\leq U^+$ always holds true, we
conclude that the two functions coincide with $U^0$, the unique
viscosity solution of \eqref{Ueffett}.

By Lemmata \ref{linkuepUep} and \ref{linkuU}, the convergence of
$U^\ep$ to $U^0$ proves in particular that $u^\ep$ converges
towards $u^0$ viscosity solution of \eqref{ueffett}.
\bigskip

To prove that  $U^+$ is a subsolution of \eqref{Ueffett}, we argue by contradiction. In what follows we will use the notation $X=(x,x_{N+1})$.  We consider a test function $\phi$ such
that $U^+-\phi$ attains a zero maximum at $(t_0,X_0)$ with $t_0>0$
and $X_0=(x_0,x_{N+1}^0)$. Without loss of generality we may
assume that the maximum is strict and global. Suppose that there
exists $\theta>0$ such that
$$\p_t\phi(t_0,X_0)=\overline{H}(\nabla_x \phi(t_0,X_0),L_0)+\theta,$$
where
\begin{equation}\label{l0}\begin{split}L_0=&\int_{|x|\leq1}(\phi(t_0,x_0+x,x_{N+1}^0)-\phi(t_0,X_0)-\nabla_x\phi(t_0,X_0)\cdot
x)\mu(dx)\\&+\int_{|x|>1}(U^+(t_0,x_0+x,x_{N+1}^0)-U^+(t_0,X_0))\mu(dx).\end{split}
\end{equation}

\noindent {\bf Step 2: Construction of $\phi^\ep$}\\
By Proposition \ref{Hprop}, we know that there exists $L_1>0$ (that we take minimal) such
that
$$\overline{H}(\nabla_x\phi(t_0,X_0),L_0)+\theta=\overline{H}(\nabla_x\phi(t_0,X_0),L_0+L_1).$$
By Propositions \ref{apprcorrectors} and \ref{Hprop}, we can
consider a sequence $L_\eta\rightarrow L_1$ as
$\eta\rightarrow0^+$, such that
$\lam^+_\eta(\nabla_x\phi(t_0,X_0),L_0+L_\eta)=\lam(\nabla_x\phi(t_0,X_0),L_0+L_1)$.
We choose $\eta$ so small that $L_\eta-o_\eta(1)\geq L_1/2>0$,
where $o_\eta(1)$ is defined in Proposition \ref{apprcorrectors}.
Let $V^+_{\eta}$ be the approximate supercorrector given by
Proposition \ref{apprcorrectors} with
$$ p=\nabla_x\phi(t_0,X_0),\quad
L=L_0+L_\eta$$ and
$$\lam^+_\eta=\lam^+_\eta(p,L_0+L_\eta)= \p_t\phi(t_0,X_0).$$
For simplicity of notations, in the following we denote
$V=V^+_\eta$. We consider the function $F(t,X)=\phi(t,X)-p\cdot
x-\lam t$, and as in \cite{im} and \cite{imr} we introduce the
"$x_{N+1}$-twisted perturbed test function" $\phi^\epsilon$
defined by:
\begin{equation}\label{phiep}\phi^\epsilon(t,X):=
\begin{cases}
\phi(t,X)+\epsilon
V\left(\frac{t}{\epsilon},\frac{x}{\epsilon},\frac{F(t,X)}{\epsilon}\right)+\epsilon
k_\epsilon
 & \text{in}\quad (\frac{t_0}{2},2t_0)\times B_\frac{1}{2}(X_0)\\
U^\epsilon (t,X) &\text{outside},
\end{cases}
\end{equation}
where $k_\epsilon\in\Z$ will be chosen later.\\

\noindent {\bf Step 3: Checking that $\phi^\ep$ is a supersolution}\\
\noindent {\bf Step 3.1: Outside $Q_{r,r}(t_0,x_0)$}\\
We are going to
prove that $\phi^\epsilon$ is a supersolution of \eqref{Uep} in
$Q_{r,r}(t_0,X_0)$ for some $r<\frac{1}{2}$ properly chosen and
such that $Q_{r,r}(t_0,X_0)\subset(\frac{t_0}{2},2t_0)\times
B_\frac{1}{2}(X_0)$. First, remark that since $U^+-\phi$ attains a
strict maximum at $(t_0,X_0)$ with $U^+-\phi=0$ at $(t_0,X_0)$ and
$V$ is bounded, we can ensure that there exists $\ep_0=\ep_0(r)>0$
such that for $\ep\leq \ep_0$
\begin{equation}\label{phiep2}U^\ep(t,X)\leq \phi(t,X)+\epsilon
V\left(\frac{t}{\epsilon},\frac{x}{\epsilon},\frac{F(t,X)}{\epsilon}\right)-\gamma_r,\quad\text{in
}\left(\frac{t_0}{3},3t_0\right)\times B_1(x_0)\setminus
Q_{r,r}(t_0,x_0)\end{equation} for some $\gamma_r=o_r(1)>0$. Hence
choosing $k_\ep=\lceil \frac{-\gamma_r}{\ep}\rceil$ we get
$U^\ep\leq \phi^\ep$ outside $Q_{r,r}(t_0,X_0)$.

\noindent {\bf Step 3.2: Inside $Q_{r_0,r_0}(t_0,x_0)$: $\phi^\ep$ tested by $\psi$}\\
Let us next study the equation. From \eqref{linklemmdis}, we
deduce that $U^+(t,x,x_{N+1}+a)=U^+(t,x,x_{N+1})+a$ for any
$a\in\R$, from  which we derive that $\p_{x_{N+1}}
F(t_0,X_0)=\p_{x_{N+1}}\phi(t_0,X_0)=1$. Then, there exists
$r_0>0$ such that the map

$$\begin{array}{cccc}
                Id\times F:&Q_{r_0,r_0}(t_0,X_0) & \longrightarrow & \mathcal{U}_{r_0} \\
                &(t,x,x_{N+1})  & \longmapsto & (t,x,F(t,x,x_{N+1})) \\
\end{array}$$
is a $C^1$-diffeomorphism from $Q_{r_0,r_0}(t_0,X_0)$ onto its
range $\mathcal{U}_{r_0}$. Let $G:\mathcal{U}_{r_0}\rightarrow\R$
be the map such that $$\begin{array}{cccc}
                Id\times G:&\mathcal{U}_{r_0} & \longrightarrow & Q_{r_0,r_0}(t_0,X_0)  \\
                &(t,x,\xi_{N+1})  & \longmapsto & (t,x,G(t,x,\xi_{N+1})) \\
\end{array}$$ is the inverse of $Id\times F$.
 Let
us introduce the variables $\tau=t/\ep$, $Y=(y,y_{N+1})$ with
$y=x/\ep$ and $y_{N+1}=F(t,X)/\ep$. Let us consider a test
function $\psi$ such that $\phi^\ep-\psi$ attains a global zero
minimum at $(\ts,\Xs)\in Q_{r_0,r_0}(t_0,X_0)$ and define \beqs
\Gamma^\ep(\tau,Y)=\frac{1}{\ep}[\psi(\ep\tau,\ep y, G(\ep\tau,\ep
y,\ep y_{N+1}))-\phi(\ep\tau,\ep y,G(\ep\tau,\ep y,\ep
y_{N+1}))]-k_\ep.\eeqs Then
$$\psi(t,X)=\phi(t,X)+\epsilon
\Gamma^\ep\left(\frac{t}{\epsilon},\frac{x}{\epsilon},\frac{F(t,X)}{\epsilon}\right)+\epsilon
k_\epsilon$$ and $\Gamma^\ep$ is a test funtion for $V$:

 \beq\label{gammaeptestv}
\Gamma^\ep(\tas,\Ys)=V(\tas,\Ys)\quad\text{and}\quad
\Gamma^\ep(\tau,Y)\leq V(\tau,Y)\quad \text{for all }(\ep\tau,\ep
Y)\in Q_{r_0,r_0}(t_0,X_0),\eeq where $\tas=\ts/\ep$,
$\ys=\xs/\ep,$ $\overline{y}_{N+1}=F(\ts,\Xs)/\ep$,
$\Ys=(\ys,\overline{y}_{N+1})$. From Proposition
\ref{apprcorrectors}, we know that $V$ is Lipschitz continuous
w.r.t. $y_{N+1}$ with Lipschitz constant $M_\eta$ depending on
$\eta$. This implies that
\beq\label{vlipyn+1}|\p_{y_{N+1}}\Gamma^\ep(\tas,\Ys)|\leq
M_\eta.\eeq

Simple computations yield with $P=(p,1)\in\R^{N+1}$:
\beq\label{equaprogfconvchangevar}\left\{%
\begin{array}{ll}
    \lam^+_\eta+\p_{\tau}\Gamma^\ep(\tas,\Ys)=\p_t\psi(\ts,\Xs)+\left(1+\p_{y_{N+1}}\Gamma^\ep(\tas,\Ys)\right)
    (\p_t\phi(t_0,X_0)-\p_t\phi(\ts,\Xs)), \\
    \lam^+_\eta \tas+P\cdot\Ys+V(\tas,\Ys)=\frac{\phi^\ep(\ts,\Xs)}{\ep}-k_\ep. \\
\end{array}%
\right.\eeq Using \eqref{equaprogfconvchangevar} and
\eqref{vlipyn+1}, Equation \eqref{apprcorrequ} yields for any
$\rho>0$
\begin{equation}\label{phiepquat}\begin{split}\p_t\psi(\ts,\Xs)+o_r(1)&\geq
L_0+L_\eta+\I^{1,\rho}[\Gamma^\ep(\tas,\cdot,\overline{y}_{N+1}),\ys]+\I^{2,\rho}[V(\tas,\cdot,\overline{y}_{N+1}),\ys]
\\&-W'\left(\frac{\phi^\ep(\ts,\Xs)}{\ep}\right)+\sigma\left(\frac{\ts}{\ep},\frac{\xs}{\ep}\right)-o_\eta(1).\end{split}\end{equation}
With the following lemma (which will be proved in the next subsection), we make rigorous the heuristic computations done in Subsection \ref{sect312} to estimate 
the error when plugging \eqref{tildeUep} in \eqref{Uepintro}.
\begin{lem} {\bf (Supersolution property for $\phi^\ep$)}\label{converglem}\\
For $\ep\leq \ep_0(r)< r\leq r_0$, we
have \beqs\begin{split}\p_t\psi(\ts,\Xs)&\geq
\I^{1,1}\left[\psi(\ts,\cdot,\overline{x}_{N+1}),\xs\right]+\I^{2,1}\left[\phi^\ep(\ts,\cdot,\overline{x}_{N+1}),\xs\right]\\&
-W'\left(\frac{\phi^\ep(\ts,\Xs)}{\ep}\right)+\sigma\left(\frac{\ts}{\ep},\frac{\xs}{\ep}\right)-o_\eta(1)+o_r(1)+L_\eta.\end{split}\eeqs\end{lem}

Let $r\leq r_0$ be so small that $o_r(1)\geq-L_1/4$. Then,
recalling that $L_\eta-o_\eta(1)\geq L_1/2$, for $\ep\leq
\ep_0(r)$ we have
\begin{equation*}\begin{split}\p_t\psi(\ts,\Xs)&\geq
\I^{1,1}\left[\psi(\ts,\cdot,\overline{x}_{N+1}),\xs\right]+\I^{2,1}\left[\phi^\ep(\ts,\cdot,\overline{x}_{N+1}),\xs\right]
-W'\left(\frac{\phi^\ep(\ts,\Xs)}{\ep}\right)\\&+\sigma\left(\frac{\ts}{\ep},\frac{\xs}{\ep}\right)+\frac{L_1}{4},\end{split}\end{equation*}
and therefore $\phi^\ep$ is a supersolution of \eqref{Uep} in
$Q_{r,r}(t_0,X_0)$.\\
 \noindent {\bf Step 4: Conclusion}\\
Since $U^\ep\leq \phi^\ep$ outside
$Q_{r,r}(t_0,X_0)$, by the comparison principle, Proposition
\ref{comparisonbounded}, we conclude that $U^\ep(t,X)\leq
\phi(t,X)+\epsilon
V\left(\frac{t}{\epsilon},\frac{x}{\epsilon},\frac{F(t,X)}{\epsilon}\right)+\ep
k_\ep$ in $Q_{r,r}(t_0,X_0)$ and we obtain the desired
contradiction by passing to the upper limit as $\ep\rightarrow 0$ at $(t_0,X_0)$ using the fact that $U^+(t_0,X_0)=\phi(t_0,X_0)$: $0\leq-\gamma_r$.\\
This ends the proof of Theorem \ref{convergence}.\\

\subsection{Proof of Lemma \ref{converglem}}
The result will follow from (\ref{phiepquat})
and the following inequality
\begin{equation}\label{eq::s15}
\begin{array}{l}
L_0+\I^{1,\rho}[\Gamma^\ep(\tas,\cdot,\overline{y}_{N+1}),\ys]+\I^{2,\rho}[V(\tas,\cdot,\overline{y}_{N+1}),\ys]\\ 
\\
\ge \I^{1,1}\left[\psi(\ts,\cdot,\overline{x}_{N+1}),\xs\right]+\I^{2,1}\left[\phi^\ep(\ts,\cdot,\overline{x}_{N+1}),\xs\right] + o_r(1)
\end{array}
\end{equation}
To show the result, we proceed in several
steps. In what follows, we denote by $C$ various positive
constants independent of $\ep$.
We start to call
$$L_0^1=\int_{|x|\leq1}(\phi(t_0,x_0+x,x^0_{N+1})-\phi(t_0,X_0)-\nabla\phi(t_0,X_0)\cdot
x)\mu(dx),$$
$$L_0^2=\int_{|x|>1}(U^+(t_0,x_0+x,x^0_{N+1})-U^+(t_0,X_0))\mu(dx).$$ Then, recalling the definition \eqref{l0} of $L_0$, we can write
\begin{equation}\begin{split}\label{L0=L1+L2}L_0=L_0^1+L_0^2.\end{split}\end{equation}
Keep in mind that $\ys_{N+1}=\frac{F(\ts,\Xs)}{\ep}$. Since
$\psi(t,X)=\phi(t,X)+\ep\Gamma^\ep\left(\frac{t}{\ep},\frac{x}{\ep},\frac{F(t,X)}{\ep}\right)+\ep
k_\ep$, we have
\begin{equation}\label{i1psi}\begin{split}\I^{1,1}\left[\psi(\ts,\cdot,\overline{x}_{N+1}),\xs\right]= I_1
+I_2,
\end{split}
\end{equation}where
$$\left\{\begin{array}{lll}
I_1 &=&\displaystyle \int_{|x|\leq1}\ep
\left(\begin{array}{l}
\Gamma^\ep\left(\frac{\overline{t}}{\ep},\frac{\xs+x}{\ep},
\frac{F(\ts,\xs+x,\overline{x}_{N+1})}{\ep}\right)-\Gamma^\ep(\tas,\Ys)\\
-\nabla_{y}\Gamma^\ep(\tas,\Ys)\cdot
\frac{x}{\ep}-\p_{y_{N+1}}\Gamma^\ep(\tas,\Ys)\nabla_x
F(\ts,\Xs)\cdot \frac{x}{\ep}
\end{array}\right)\mu(dx),\\
&&\\
I_2& =&\displaystyle \int_{|x|\leq 1}\left(\phi(\ts,\xs+x,\overline{x}_{N+1})-\phi(\ts,\Xs)-\nabla\phi(\ts,\Xs)\cdot
x\right)\mu(dx).
\end{array}\right.$$
In order to show (\ref{eq::s15}), we show successively in Steps 1, 2 and 3:
$$\left\{\begin{array}{l}
I_1\leq
\I^{1,\rho}[\Gamma^\ep(\tas,\cdot,\overline{y}_{N+1}),\ys]+\I^{2,\rho}[V(\tas,\cdot,\overline{y}_{N+1}),\ys]
+o_r(1)+C_\ep\rho\\
\\
I_2\leq L_0^1+ o_r(1)\\
\\
\I^{2,1}\left[\phi^\ep(\ts,\cdot,\overline{x}_{N+1}),\xs\right]\leq L_0^2+ o_r(1)
\end{array}\right.$$
Because the expressions are non linear and non local and with a singular kernel, 
there is no simple computation and we have to carefully check those inequalities 
sometimes splitting terms in easier parts to estimate.

\noindent{\bf Step 1:} We can choose $\ep_0$ so small that for any
$\ep\leq\ep_0$ and any $\rho>0$ small enough $$I_1\leq
\I^{1,\rho}[\Gamma^\ep(\tas,\cdot,\overline{y}_{N+1}),\ys]+\I^{2,\rho}[V(\tas,\cdot,\overline{y}_{N+1}),\ys]
+o_r(1)+C_\ep\rho.$$ \\

Take $\rho>0$,  $\delta>\rho$ small and $R>0$ large and such that
$\ep R<1$. Since $g$ is even, we can write

$$I_1=I_1^0+I_1^1+I_1^2+I_1^3,$$ where
\begin{equation*}\begin{split}I_1^0&=\int_{|x|\leq\ep\rho}\ep\left(\Gamma^\ep\left(\frac{\overline{t}}{\ep},\frac{\xs+x}{\ep},
\frac{F(\ts,\xs+x,\overline{x}_{N+1})}{\ep}\right)-\Gamma^\ep(\tas,\Ys)-\nabla_{y}\Gamma^\ep(\tas,\Ys)\cdot
\frac{x}{\ep}\right.\\&\left.-\p_{y_{N+1}}\Gamma^\ep(\tas,\Ys)\nabla_x
F(\ts,\Xs)\cdot \frac{x}{\ep}\right)\mu(dx),\end{split}
\end{equation*}
$$I_1^1=\int_{\ep\rho\leq|x|\leq\ep \delta}\ep\left(\Gamma^\ep\left(\frac{\overline{t}}{\ep},\frac{\xs+x}{\ep},
\frac{F(\ts,\xs+x,\overline{x}_{N+1})}{\ep}\right)-
\Gamma^\ep(\tas,\Ys)\right)\mu(dx),$$
$$I_1^2=\int_{\ep\delta\leq|x|\leq\ep R}\ep\left(\Gamma^\ep\left(\frac{\overline{t}}{\ep},\frac{\xs+x}{\ep},
\frac{F(\ts,\xs+x,\overline{x}_{N+1})}{\ep}\right)-
\Gamma^\ep(\tas,\Ys)\right)\mu(dx),$$
$$I_1^3=\int_{\ep R\leq|x|\leq1}\ep\left(\Gamma^\ep\left(\frac{\overline{t}}{\ep},\frac{\xs+x}{\ep},
\frac{F(\ts,\xs+x,\overline{x}_{N+1})}{\ep}\right)-
\Gamma^\ep(\tas,\Ys)\right)\mu(dx).$$

 Moreover
$$\I^{2,\rho}[V(\tas,\cdot,\overline{y}_{N+1}),\ys]=J_1+J_2+J_3,$$ where

$$J_1=\int_{\rho<|z|\leq\delta}(V(\tas,\ys+z,\overline{y}_{N+1})- V(\tas,\Ys))\mu(dz),$$
$$J_2=\int_{\delta<|z|\leq R}(V(\tas,\ys+z,\overline{y}_{N+1})- V(\tas,\Ys))\mu(dz),$$

$$J_3=\int_{|z|>R}(V(\tas,\ys+z,\overline{y}_{N+1})- V(\tas,\Ys))\mu(dz).$$

\noindent STEP 1.1: \emph{Estimate of $I_1^0$ and
$\I^{1,\rho}[\Gamma^\ep(\tas,\cdot,\overline{y}_{N+1}),\ys]$.}

Since $\Gamma^\ep$ is of class $C^2$, we have \beq\label{step1.1}
|I_1^0|,\,|\I^{1,\rho}[\Gamma^\ep(\tas,\cdot,\overline{y}_{N+1}),\ys]|\leq
C_\ep\rho,\eeq where $C_\ep$ depends on the second derivatives of
$\Gamma^\ep$. Remark that if we knew that $V$ is smooth in $y$ too, we could choose $\rho=0$. 

\noindent STEP 1.2{\emph{ Estimate of $I_1^1-J_1$. }

Using \eqref{gammaeptestv} and the fact that $g$ is even, we can
estimate $I_1^1-J_1$ as follows
\begin{equation*}\begin{split}I_1^1-J_1&\leq\int_{\rho<|z|\leq\delta}\left[V\left(\tas,\ys+z,\frac{F(\ts,\xs+\ep
z,\xs_{N+1})}{\ep}\right)-V\left(\tas,\ys+z,\frac{F(\ts,\xs)}{\ep}\right)\right]\mu(dz)\\&
=\int_{\rho<|z|\leq\delta}\left\{\left[V\left(\tas,\ys+z,\frac{F(\ts,\xs+\ep
z,\xs_{N+1})}{\ep}\right)-V\left(\tas,\ys+z,\frac{F(\ts,\xs)}{\ep}\right)\right.\right.\\&\left.
-\p_{y_{N+1}}V\left(\tas,\ys+z,\frac{F(\ts,\Xs)}{\ep}\right)\nabla_xF(\ts,\Xs)\cdot
z\right]\\&
\left.+\left[\p_{y_{N+1}}V(\tas,\ys+z,\ys_{N+1})-\p_{y_{N+1}}V(\tas,\Ys)\right]\nabla_x
F(\ts,\Xs)\cdot z\right\}\mu(dz).
\end{split}
\end{equation*}
Next, using \eqref{appcorr3} and \eqref{contrderivappcorr}, we get
\beq\label{i1-j1}I_1^1-J_1\leq
C\int_{|z|\leq\delta}(|z|^2+|z|^{1+\al})\mu(dz)\leq C
\delta^\al.\eeq

\noindent STEP 1.3{\emph{ Estimate of $I_1^2-J_2$. }

If $M_\eta$ is the Lipschitz constant of $V$ w.r.t. $y_{N+1}$,
then
\begin{equation*}\begin{split}I_1^2-J_2&\leq\int_{\delta<|z|\leq
R}\left(V\left(\tas,\ys+z,\frac{F(\ts,\xs+\ep
z,\xs_{N+1})}{\ep}\right)-V\left(\tas,\ys+z,\frac{F(\ts,\Xs)}{\ep}\right)\right)\mu(dz)\\
&\leq M_\eta\int_{\delta<|z|\leq R}\left|\frac{F(\ts,\xs+\ep
z,\xs_{N+1})}{\ep}-\frac{F(\ts,\Xs)}{\ep}\right|\mu(dz)\\&\leq
M_\eta\int_{\delta<|z|\leq R}\sup_{|z|\leq
R}|\nabla_xF(\ts,\xs+\ep z,\xs_{N+1})||z|\mu(dz).
\end{split}
\end{equation*}
Then \beq\label{i2-j2}I_1^2-J_2\leq C\sup_{|z|\leq
R}|\nabla_xF(\ts,\xs+\ep z,\xs_{N+1})|\log(R/\delta)\eeq

\noindent STEP 1.4: {\em Estimate of $I_1^3$ and $J_3$. }

Since $V$ is uniformly bounded on $\R^+\times\R^{N+1}$, we have
\begin{equation}\label{i3est}\begin{split}I_1^3&\leq
\int_{R<|z|\leq
\frac{1}{\ep}}\left(V\left(\tas,\ys+z,\frac{F(\ts,\xs+\ep
z,\xs_{N+1})}{\ep}\right)- V(\tas,\Ys)\right)\mu(dz)\\&\leq
\int_{|z|>R}2\|v\|_\infty\mu(dz)\leq \frac{C}{R}.
\end{split}
\end{equation}
Similarly \beq\label{jeest}|J_3|\leq \frac{C}{R}.\eeq

Now, from \eqref{step1.1}, \eqref{i1-j1}, \eqref{i2-j2},
\eqref{i3est} and \eqref{jeest}, we infer that \beqs\begin{split}
I_1&\leq
\I^{1,\rho}[\Gamma^\ep(\tas,\cdot,\overline{y}_{N+1}),\ys]+\I^{2,\rho}[V(\tas,\cdot,\overline{y}_{N+1}),\ys]+
2C_\ep\rho+C \delta^\al\\&+C\sup_{|z|\leq R}|\nabla_xF(\ts,\xs+\ep
z,\xs_{N+1})|\log\left(\frac{R}{\delta}\right)+\frac{C}{R}.\end{split}\eeqs

We choose $R=R(r)$ such $R\rightarrow+\infty$ as
$r\rightarrow0^+$, $\ep_0=\ep_0(r)$ such that $R\ep_0(r)\leq r$
and $\delta=\delta(r)>0$ such that $\delta\rightarrow0$ as
$r\rightarrow0^+$ and $r\log(R/\delta)\rightarrow0$ as
$r\rightarrow0^+$. With this choice, for any $\ep\leq\ep_0$ and
any $\rho<\delta$
$$C \delta^\al+C\sup_{|z|\leq
R}|\nabla_xF(\ts,\xs+\ep
z,\xs_{N+1})|\log\left(\frac{R}{\delta}\right)+\frac{C}{R}=o_r(1)\quad\text{as
}r\rightarrow0^+,$$ and Step 1 is proved.

\noindent{\bf Step 2: $I_2\leq L_0^1+ o_r(1)$. }

For $0<\nu<1$ we can split $I_2$ and $L_0^1$ as follows
\beqs\begin{split}I_2&=\int_{|x|\leq\nu}(\phi(\ts,\xs+x,\xs_{N+1})-\phi(\ts,\Xs)-\nabla\phi(\ts,\Xs)\cdot
x)\mu(dx)\\&+\int_{\nu\leq|x|\leq1}(\phi(\ts,\xs+x,\xs_{N+1})-\phi(\ts,\Xs))\mu(dx)=I_2^1+I_2^2
,\end{split}\eeqs
\beqs\begin{split}L_0^1&=\int_{|x|\leq\nu}(\phi(t_0,x_0+x,x^0_{N+1})-\phi(t_0,X_0)-\nabla\phi(t_0,X_0)\cdot
x)\mu(dx)\\&+\int_{\nu\leq|x|\leq1}(\phi(t_0,x_0+x,x^0_{N+1})-\phi(t_0,X_0))\mu(dx)
=T_1+T_2.\end{split}\eeqs Since $\phi$ is of class $C^2$ we have
$$I_2^1,T_1\leq C\nu.$$Using the Lipschitz continuity of $\phi$ we
get
\begin{equation*}\begin{split}I_2^2-T_2=\int_{\nu<|x|\leq1}Cr\mu(dx)\leq
C\frac{r}{\nu}.\end{split}\end{equation*}
 Hence, Step 2 follows choosing $\nu=\nu(r)$ such that
 $\nu\rightarrow0$ and $r/\nu\rightarrow0$ as $r\rightarrow0^+$.

 \noindent{\bf Step 3: $\I^{2,1}\left[\phi^\ep(\ts,\cdot,\overline{x}_{N+1}),\xs\right]\leq L_0^2+ o_r(1)$. }

Remark that
$$U^\ep(\ts,\xs+x,\xs_{N+1})-\phi(\ts,\Xs)-\ep V(\tas,\Ys)-\ep k_\ep\leq
U^+(t_0,x_0+x,x^0_{N+1})-\phi(t_0,X_0)+o_\ep(1)+o_r(1).$$Then,
recalling that $\phi(t_0,X_0)=U^+(t_0,X_0),$ for $\ep\leq\ep_0$ we
get
$$\I^{2,1}\left[\phi^\ep(\ts,\cdot,\overline{x}_{N+1}),\xs\right]-L_0^2\leq o_r(1)$$
and Step 3 is proved.

Finally  \eqref{L0=L1+L2}, \eqref{i1psi}, Steps 1, 2 and 3 give
\beqs\begin{split}
\I^{1,1}\left[\psi(\ts,\cdot,\overline{x}_{N+1}),\xs\right]+\I^{2,1}\left[\phi^\ep(\ts,\cdot,\overline{x}_{N+1}),\xs\right]&\leq
\I^{1,\rho}[\Gamma^\ep(\tas,\cdot,\overline{y}_{N+1}),\ys]+\I^{2,\rho}[V(\tas,\cdot,\overline{y}_{N+1}),\ys]\\&+L_0
+o_r(1)+C_\ep\rho.\end{split}\end{equation*} from which, using
inequality \eqref{phiepquat} and letting $\rho\rightarrow0^+$, we
get for $\ep\leq\ep_0$ \beqs\begin{split}\p_t\psi(\ts,\Xs)&\geq
\I^{1,1}\left[\psi(\ts,\cdot,\overline{x}_{N+1}),\xs\right]+\I^{2,1}\left[\phi^\ep(\ts,\cdot,\overline{x}_{N+1}),\xs\right]
-W'\left(\frac{\phi^\ep(\ts,\Xs)}{\ep}\right)+\sigma\left(\frac{\ts}{\ep},\frac{\xs}{\ep}\right)\\&-o_\eta(1)+o_r(1)+L_\eta\end{split}\end{equation*}
and this concludes the proof of the lemma. \finedim

\section{Building of Lipschitz sub and supercorrectors}\label{lipcorrsec}

In this section we construct  bounded sub and supersolutions of \eqref{V}
that are Lipschitz w.r.t. $y_{N+1}$. As a byproduct, we will prove
Theorem \ref{ergodic} and Proposition \ref{Hprop}.
\begin{prop}[Lipschitz continuous sub and supercorrectors]\label{lipcorrect}Let $\lam$ be the quantity defined by Theorem \ref{ergodic}. Then,
for any fixed $p\in\R^N$, $P=(p,1)$, $L\in\R$ and $\eta>0$ small
enough, there exist real numbers $\lam^+_\eta(p,L)$,
$\lam^-_\eta(p,L)$, a constant $C>0$ (independent of $\eta,\,p$
and $L$) and bounded super and subcorrectors $W^+_{\eta},
W^-_{\eta}$ i.e. respectively a super and a subsolution of
\eqref{V} (with respectively $\lam^+_\eta$ and $\lam^-_\eta$ in
place of $\lam$) such that \beqs
\lim_{\eta\rightarrow0^+}\lam^+_\eta(p,L)=\lim_{\eta\rightarrow0^+}\lam^-_\eta(p,L)=\lam(p,L),\eeqs
$\lam^{\pm}_\eta$ satisfy (i) and (ii) of Proposition \ref{Hprop}
and for any $(\tau,Y)\in\R^+\times\R^{N+1}$
\beq\label{lipcorrbound}|W^{\pm}_{\eta}(\tau,Y)|\leq C.\eeq
 Moreover $W^{\pm}_{\eta}$ are Lipschitz continuous w.r.t. $y_{N+1}$ and $\al$-H\"{o}lder
continuous w.r.t. $y$ for any $0<\al<1$, with  \beq\label{lipyn+1}
-1\leq\p_{y_{N+1}}W^{\pm}_\eta \leq
\frac{\|W''\|_\infty}{\eta},\eeq \beq\label{contrderivappcorrlip}
<W^{\pm}_\eta>_y^\al\leq C_{\eta}.\eeq
\end{prop}
 In order to prove the proposition, for $\eta\geq0$, $a_0,\,L\in\R$,
$p\in\R^N$ and $P=(p,1)$, we introduce the problem
\begin{equation}\label{wlip}\left\{
  \begin{array}{ll}
    \p_{\tau} U=L+\I[U(\tau,\cdot,y_{N+1})]-W'(U+P\cdot Y)+\s(\tau,y)\\
    \quad\quad\,+\eta[a_0+\inf_{Y'}U(\tau,Y')-U(\tau,Y)]|\p_{y_{N+1}}U+1|
    & \hbox{in } \R^+\times\R^{N+1}\\
    U(0,Y)=0 & \hbox{on }\R^{N+1}.
  \end{array}
\right.\end{equation}

 We have the following result whose proof is postponed to the Appendix (Section \ref{appendix}).
\begin{prop}[Comparison principle for \eqref{wlip}]\label{pro::s19}
Let $U_1\in USC_b(\R^+\times\R^{N+1})$ and $U_2\in LSC_b(\R^+\times\R^{N+1})$ be respectively
a viscosity subsolution and supersolution of \eqref{wlip}, then
$U_1\leq U_2$ on $\R^+\times\R^{N+1}$.
\end{prop}

\subsection{Lipschitz regularity}

\begin{prop}[Lipschitz continuity in $y_{N+1}$] Suppose $\eta>0$. Let $U_\eta\in C_b(\R^+\times\R^{N+1})$ be the viscosity solution of \eqref{wlip}.
Then $U_\eta$ is Lipschitz continuous w.r.t. $y_{N+1}$ and for
almost every $(\tau,Y)\in \R^+\times\R^{N+1}$
 \beq\label{derivUyn+1}-1\leq \p_{y_{N+1}}U_\eta(\tau,Y)\leq
\frac{\|W''\|_\infty}{\eta}.\eeq
\end{prop}

For a formal argument, we refer the reader to Step 1 of  Subsection \ref{sect3.2}.\\

\dim Let us define $\widehat{U}(\tau,Y)=U(\tau,Y)+y_{N+1}$, then
$\widehat{U}$ satisfies
\begin{equation}\label{ulip}\left\{
  \begin{array}{ll}
    \p_{\tau} \widehat{U}=L+\I[\widehat{U}(\tau,\cdot,y_{N+1})]-W'(\widehat{U}+p\cdot y)+\s(\tau,y)\\
    \quad\quad\,+\eta[a_0+\inf_{Y'}(\widehat{U}(\tau,Y')-y'_{N+1})-(\widehat{U}(\tau,Y)-y_{N+1})]|\p_{y_{N+1}}\widehat{U}|
    & \hbox{in } \R^+\times\R^{N+1}\\
    \widehat{U}(0,Y)=y_{N+1} & \hbox{on }\R^{N+1}.
  \end{array}
\right.\end{equation} We are going to prove that $\widehat{U}$ is
Lipschitz continuous w.r.t. $y_{N+1}$ with \beqs 0\leq
\p_{y_{N+1}}\widehat{U}(\tau,Y)\leq
1+\frac{\|W''\|_\infty}{\eta}.\eeqs By comparison,
$\widehat{U}(t,y,y_{N+1})\leq \widehat{U}(t,y,y_{N+1}+h)$ for
$h\geq0$, from which immediately follows that
$\p_{y_{N+1}}\widehat{U}\geq 0$. In particular we can replace
$|\p_{y_{N+1}}\widehat{U}|$ by $\p_{y_{N+1}}\widehat{U}$ in
\eqref{ulip}.

 Let us now show that $\p_{y_{N+1}}\widehat{U}\leq
1+\frac{\|W''\|_\infty}{\eta}$. We argue by contradiction by
assuming that for some $T>0$ the supremum of the function
$\widehat{U}(\tau,y,y_{N+1})-\widehat{U}(\tau,y,z_{N+1})-K|\yn-\zn|$
on $[0,T]\times\R^{N+1}$ is strictly positive as soon as
$K>1+\frac{\|W''\|_\infty}{\eta}$. Then for $\delta,\beta>0$ small
enough, $M$ defined by
$$M=\max_{(\tau,y)\in[0,T]\times\R^N\atop y_{N+1},z_{N+1}\in\R }
\left(\widehat{U}(\tau,y,y_{N+1})-\widehat{U}(\tau,y,z_{N+1})-K|y_{N+1}-z_{N+1}|-\beta\psi(Y)-\frac{\delta}{T-\tau}\right),$$
where $\psi$ is defined as the function $\psi_2$ in the proof of
Proposition \ref{regularityvisc}, is positive. For $j>0$ let
\beqs\begin{split}M_j&=\max_{\tau,s\in[0,T],y,z\in\R^N\atop
y_{N+1},z_{N+1}\in\R }
\left(\widehat{U}(\tau,y,y_{N+1})-\widehat{U}(s,z,z_{N+1})-K|y_{N+1}-z_{N+1}|-\beta\psi(Y)\right.\\&\left.-\frac{\delta}{T-\tau}-j|\tau-s|^2
-j|y-z|^2\right),\end{split}\eeqs and let
$(\tau^j,y^j,y_{N+1}^j,s^j,z^j,z_{N+1}^j)\in
([0,T]\times\R^{N+1})^2$ be a point where $M_j$ is attained.
Classical arguments show that $M_j\rightarrow M$,
$(\tau^j,y^j,y_{N+1}^j,s^j,z^j,z_{N+1}^j)\rightarrow(\tas,
\ys,\overline{y}_{N+1},\tas,\ys,\overline{z}_{N+1})$ as
$j\rightarrow+\infty$, where $(\tas,\ys,
\overline{y}_{N+1},\overline{z}_{N+1})$ is a point where $M$ is
attained.

Remark that $0<\tas<T$, moreover, since
$\widehat{U}(\tas,\ys,\overline{y}_{N+1})>\widehat{U}(\tas,\ys,\overline{z}_{N+1})$
and  $\widehat{U}$ is nondecreasing in $y_{N+1}$, it is
\beq\label{yn>zn}\overline{y}_{N+1}>\overline{z}_{N+1}.\eeq In
particular $y_{N+1}^j\neq z_{N+1}^j$ and $0<s_j,\,\tau_j<T$ for
$j$ large enough. Hence, for $r>0$, we obtain the following
viscosity inequalities
\beq\label{5ulipusub}\begin{split}&\frac{\delta}{(T-\tau_j)^2}+j(t_j-s_j)\\&\leq
L+C_Njr+\beta\I^{1,r}[\psi(\cdot,y_{N+1}^j),y^j]+\I^{2,r}[\widehat{U}(\tau^j,\cdot,y_{N+1}^j),y^j]\\&-W'(\widehat{U}(\tau^j,y^j,y_{N+1}^j)+p\cdot
y^j)+\s(\tau^j,y^j)+\eta[a_0+\inf_{Y'}(\widehat{U}(\tau_j,Y')-y'_{N+1})\\&-(\widehat{U}(\tau^j,y^j,y_{N+1}^j)-y_{N+1}^j)]\left(K\frac{y_{N+1}^j-z_{N+1}^j
}{|y_{N+1}^j-z_{N+1}^j|}+\beta\p_{y_{N+1}}\psi(y^j,y_{N+1}^j)\right)
,\end{split}\eeq and
\beq\label{5ulipusuper}\begin{split}&j(t_j-s_j)\\&\geq
L-C_Njr+\I^{2,r}[\widehat{U}(s^j,\cdot,z_{N+1}^j),z^j]-W'(\widehat{U}(s^j,z^j,z_{N+1}^j)+p\cdot
z^j)+\s(s^j,z^j)\\&+\eta[a_0+\inf_{Y'}(\widehat{U}(s_j,Y')-y'_{N+1})-(\widehat{U}(s^j,z^j,z_{N+1}^j)-z_{N+1}^j)]K\frac{y_{N+1}^j-z_{N+1}^j
}{|y_{N+1}^j-z_{N+1}^j|},\end{split}\eeq where $C_N$ is a constant
depending on $N$. Since $(\tau^j,y^j,y_{N+1}^j,s^j,z^j,z_{N+1}^j)$
is a maximum point,  we have
\beqs\begin{split}\widehat{U}(\tau^j,y^j+x,y_{N+1}^j)-\widehat{U}(\tau^j,y^j,y_{N+1}^j)&\leq
\widehat{U}(s^j,z^j+x,z_{N+1}^j)-\widehat{U}(s^j,z^j,z_{N+1}^j)\\&+\beta[\psi(y^j+x,y_{N+1}^j)-\psi(y^j,y_{N+1}^j)]\end{split}\eeqs
for any $x\in\R^N$, which implies that for any $r>0$
$$\I^{2,r}[\widehat{U}(\tau^j,\cdot,y_{N+1}^j),y^j]\leq
\I^{2,r}[\widehat{U}(s^j,\cdot,z_{N+1}^j),z^j]+\beta\I^{2,r}[\psi(\cdot,y_{N+1}^j),y^j].$$
Hence, subtracting \eqref{5ulipusub} with \eqref{5ulipusuper},
sending $r\rightarrow0^+$ and then $j\rightarrow+\infty$, we get
\beqs\begin{split}\frac{\delta}{(T-\tas)^2}&\leq
\beta\I[\psi(\cdot,\overline{y}_{N+1}),\ys]+
W'(\widehat{U}(\tas,\ys,\overline{z}_{N+1})+p\cdot
\ys)-W'(\widehat{U}(\tas,\ys,\overline{y}_{N+1})+p\cdot
\ys)\\&-\eta[\widehat{U}(\tas,\ys,\overline{y}_{N+1})-\widehat{U}(\tas,\ys,\overline{z}_{N+1})-(\overline{y}_{N+1}-\overline{z}_{N+1})]
K\frac{\overline{y}_{N+1}-\overline{z}_{N+1}}{|\overline{y}_{N+1}-\overline{z}_{N+1}|}\\&
+\beta\p_{y_{N+1}}\psi(\ys,\overline{y}_{N+1})\eta[a_0+\inf_{Y'}(\widehat{U}(\tas,Y')-y'_{N+1})
-(\widehat{U}(\tas,\ys,\overline{y}_{N+1})-\overline{y}_{N+1})]
\\&\leq
\|W''\|_\infty|\widehat{U}(\tas,\ys,\overline{y}_{N+1})-\widehat{U}(\tas,\ys,\overline{z}_{N+1})|
\\&-K\eta[\widehat{U}(\tas,\ys,\overline{y}_{N+1})-\widehat{U}(\tas,\ys,\overline{z}_{N+1})-(\overline{y}_{N+1}-\overline{z}_{N+1})]
\frac{\overline{y}_{N+1}-\overline{z}_{N+1}}{|\overline{y}_{N+1}-\overline{z}_{N+1}|}
+\beta C.
\end{split}\eeqs
Then, using \eqref{yn>zn} and that
$K|\overline{y}_{N+1}-\overline{z}_{N+1}|<\widehat{U}(\tas,\ys,\overline{y}_{N+1})-\widehat{U}(\tas,\ys,\overline{z}_{N+1})$,
for $\beta$ small enough, we finally obtain

$$(\|W''\|_\infty+\eta-\eta
K)(\widehat{U}(\tas,\ys,\overline{y}_{N+1})-\widehat{U}(\tas,\ys,\overline{z}_{N+1}))\geq0,$$
which is a contradiction for $K>1+\frac{\|W''\|_\infty}{\eta}$.
\finedim

\subsection{Ergodicity}\label{ergodicitysec}
\begin{prop}[Ergodic properties]\label{ergodic2}
There exists a unique $\lam_\eta=\lam_\eta(p,L)$ such that the
viscosity solution $U_\eta\in C_b(\R^+\times\R^{N+1})$  of
\eqref{wlip}  with $\eta\geq 0$, satisfies:
\begin{equation}\label{w-lam}|U_\eta(\tau,Y)-\lam_\eta\tau|\leq
C_3\text{ for all }\tau>0,\,Y\in \R^{N+1},\end{equation}with $C_3$
independent of $\eta$.
Moreover
\begin{equation}\label{lambounds}L-\|W'\|_\infty-\|\s\|_\infty+\eta a_0\leq\lam_\eta\leq
L+\|W'\|_\infty+\|\s\|_\infty+\eta a_0.\end{equation}
\end{prop}
\dim For simplicity of notations, in what follows we denote
$U=U_\eta$ and $\lam=\lam_\eta$.

To prove the proposition we follow the proof of the analogue
result in \cite{imr}. We
proceed in three steps.\\

\noindent{\bf Step 1: existence} The functions
$W^+(\tau,Y)=C^+\tau$ and $W^-(\tau,Y)=C^-\tau$, where
$$C^{\pm}=L\pm\|W'\|_\infty\pm\|\s\|_\infty+\eta a_0,$$ are respectively sub and supersolution of \eqref{wlip}. Then the existence
of a unique solution of \eqref{wlip} follows from Perron's method.

 \noindent{\bf Step 2: control of the oscillations w.r.t. space.}\\We want
 to prove that there exists $C_1>0$ such that
 \begin{equation}\label{w3}|U(\tau,Y)-U(\tau,Z)|\leq
C_1\quad\text{for all }\tau\geq 0,\, Y,Z\in\R^{N+1}.\end{equation}
STEP 2.1. 
For a given $k\in \Z^{N+1}$, we set $P\cdot k=l+\al$, with
$l\in\Z$ and $\al\in[0,1)$. The function
$\widetilde{U}(\tau,Y)=U(\tau,Y+k)+\al$ is still a solution of
\eqref{wlip}, with $\widetilde{U}(0,Y)=\al$
 Moreover
$$U(0,Y)=0\leq \widetilde{U}(0,Y)=\al\leq 1=U(0,Y)+1.$$
Then from the comparison principle for \eqref{wlip} and invariance
by integer translations we deduce for all $\tau\geq 0$:
\begin{equation}\label{w2}|U(\tau,Y+k)-U(\tau,Y)|\leq 1.\end{equation}

\noindent STEP 2.2. We proceed as in \cite{imr} by considering the
functions
$$M(\tau):=\sup_{Y\in\R^{N+1}}U(\tau,Y),\quad
m(\tau):=\inf_{Y\in\R^{N+1}}U(\tau,Y),$$
$$q(\tau):=M(\tau)-m(\tau)=\text{osc }U(\tau,\cdot).$$

Let us assume that the extrema defining these functions are
attained: $M(\tau)=U(\tau,Y^{\tau})$, $m(\tau)=U(\tau,Z^{\tau})$.

It is easy to see that $M(\tau)$ and $m(\tau)$ satisfy in the
viscosity sense
\begin{equation*}  \p_\tau M \leq L+\I^{2}[U(\tau,\cdot,y_{N+1}^\tau),y^\tau]-W'(M+P\cdot
Y^\tau)+\s(\tau,y^\tau)+\eta[a_0+m(\tau)-M(\tau)],\end{equation*}
\begin{equation*}\p_\tau m\geq L+
\I^{2}[U(\tau,\cdot,z_{N+1}^\tau),z^\tau]-W'(m+P\cdot
Z^\tau)+\s(\tau,z^\tau)+\eta a_0.\end{equation*}

 Then $q$ satisfies in the viscosity
sense
\begin{equation*}\begin{split}\p_\tau q&\leq
\I^{2}[U(\tau,\cdot,y_{N+1}^\tau),y^\tau]-\I^{2}[U(\tau,\cdot,z_{N+1}^\tau),z^\tau]-W'(M+P\cdot
Y^\tau)\\&+W'(m+P\cdot Z^\tau)+\s(\tau,y^\tau)
-\s(\tau,z^\tau)\\&\leq
\I^{2}[U(\tau,\cdot,y_{N+1}^\tau),y^\tau]-\I^{2}[U(\tau,\cdot,z_{N+1}^\tau),z^\tau]+2\|W'\|_\infty+2\|\s\|_\infty.\end{split}\end{equation*}
Let us estimate the quantity
$\mathcal{L}(\tau):=\I^{2}[U(\tau,\cdot,y_{N+1}^\tau),y^\tau]-\I^{2}[U(\tau,\cdot,z_{N+1}^\tau),z^\tau]$
from above by a function of $q$. Let us define $k^\tau\in\Z^{N+1}$
such that $Y^\tau-(Z^\tau+k^\tau)\in[0,1)^{N+1}$ and let
$\widetilde{Z}^\tau:=Z^\tau+k^\tau$. Using successively \eqref{w2}
and the first inequality in \eqref{derivUyn+1}, we obtain:
\begin{equation*}\begin{split}\mathcal{L}(\tau)&
\leq
\int_{|z|>1}(U(\tau,y^\tau+z,y^\tau_{N+1})-U(\tau,Y^\tau))\mu(dz)
\\&-\int_{|z|>1}(U(\tau,\widetilde{z}^\tau+z,\widetilde{z}^\tau_{N+1})-U(\tau,Z^\tau))\mu(dz)+\overline{\mu}\\&
\leq
\int_{|z|>1}(U(\tau,y^\tau+z,y^\tau_{N+1})-U(\tau,Y^\tau))\mu(dz)
\\&-\int_{|z|>1}(U(\tau,\widetilde{z}^\tau+z,y^\tau_{N+1})-U(\tau,Z^\tau))\mu(dz)+2\overline{\mu},\end{split}\end{equation*}
where $\overline{\mu}=\|\mu_0\|_{L^1(\R^N\setminus B_1(0))}.$ Now,
let us introduce $c^\tau=\frac{y^\tau+\widetilde{z}^\tau}{2}$ and
$\delta^\tau=\frac{y^\tau-\widetilde{z}^\tau}{2}\in
[0,\frac{1}{2})^N$ so that $y^\tau=c^\tau+\delta^\tau$ and
$\widetilde{z}^\tau=c^\tau-\delta^\tau$. Hence
\begin{equation*}\begin{split}\mathcal{L}(\tau)&\leq
2\overline{\mu}+\int_{|z|>1}(U(\tau,c^\tau+z+\delta^\tau,y^\tau_{N+1})-U(\tau,Y^\tau))\mu(dz)
\\&-\int_{|z|>1}(U(\tau,c^\tau+z-\delta^\tau,y^\tau_{N+1})-U(\tau,Z^\tau))\mu(dz)
\\&
\leq2\overline{\mu}+\int_{|z-\delta^\tau|>1}(U(\tau,c^\tau+z,y^\tau_{N+1})-U(\tau,Y^\tau))
\mu_0(z-\delta^\tau)dz
\\&-\int_{|z+\delta^\tau|>1}(U(\tau,c^\tau+z,y^\tau_{N+1})-U(\tau,Z^\tau))\mu_0(z+\delta^\tau)
dz\\& \leq
2\overline{\mu}-\int_{\{|z-\delta^\tau|>1\}\cap\{|z+\delta^\tau|>1\}}(U(\tau,Y^\tau)-U(\tau,Z^\tau))
\min\{\mu_0(z-\delta^\tau),\mu_0(z+\delta^\tau)\}dz\\&\leq
2\overline{\mu}-c_0q(\tau)\end{split}\end{equation*} where
$c_0>0$. We conclude that $q$ satisfies in the viscosity
sense$$\p_\tau q(\tau)\leq
2\|W'\|_\infty+2\|\s\|_\infty+2\overline{\mu}-c_0q(\tau),$$ with
$q(0)=0$, from which we obtain  \eqref{w3}.

If the extrema are not attained, it suffices to consider for
$\beta>0$,
$M_\beta(\tau):=\sup_{Y\in\R^{N+1}}(U(\tau,Y)-\beta\psi(Y))$,
$m_\beta(\tau):=\inf_{Y\in\R^{N+1}}(U(\tau,Y)+\beta\psi(Y)),$ and
$q_\beta(\tau):=M_\beta(\tau)-m_\beta(\tau)$, where $\psi$ is
defined as the function $\psi_2$ in the proof of Proposition
\ref{regularityvisc}. By the properties of $\psi$, $M_\beta(\tau)$
and $m_\beta(\tau)$ are attained. Then, the previous argument
shows that
$$q_\beta\leq C_1+C\beta,$$ and passing to the limit as
$\beta\rightarrow0^+$ we get \eqref{w3}.
\\

{\bf Step 3: control of the oscillations in time.} We follow
\cite{imr} by introducing the two quantities:
$$\lam^+(T):=\sup_{\tau\geq0}\frac{U(\tau+T,0)-U(\tau,0)}{T}\quad\text{and}\quad\lam^-(T):=\inf_{\tau\geq0}\frac{U(\tau+T,0)-U(\tau,0)}{T},$$
and proving that they have a common limit as
$T\rightarrow+\infty$. First let us estimate $\lam^+(T)$ from
above. The function $U^+(t,Y):=U(\tau,0)+C_1+C^+t$, is a
supersolution of \eqref{wlip} if
$C^+=L+\|W'\|_\infty+\|\s\|_\infty+\eta a_0$. Since $U^+(0,Y)\geq
U(\tau,Y)$ if $C_1$ is as in \eqref{w3}, by the comparison
principle for \eqref{wlip} in the time interval
$[\tau,\tau+\tau_0]$, for any $\tau_0>0$ and $t\in[0,\tau_0]$ we
get
\begin{equation}\label{w4}U(\tau+t,Y)\leq
U(\tau,0)+C_1+C^+t.\end{equation} Similarly
\begin{equation}\label{w5}U(\tau+t,Y)\geq
U(\tau,0)-C_1+C^-t,\end{equation}where
$C^-=L-\|W'\|_\infty-\|\s\|_\infty+\eta a_0$. We then obtain for
$\tau_0=t=T$ and $y=0$:
\begin{equation}\label{lam+-}L-\|W'\|_\infty-\|\s\|_\infty+\eta a_0-\frac{C_1}{T}\leq\lam^-(T)\leq\lam^+(T)\leq
L+\|W'\|_\infty+\|\s\|_\infty+\eta
a_0+\frac{C_1}{T}.\end{equation} By definition of $\lam^{\pm}(T)$,
for any $\delta>0$, there exist $\tau^{\pm}\geq0$ such that
$$\left|\lam^{\pm}(T)-\frac{U(\tau^{\pm}+T,0)-U(\tau^{\pm},0)}{T}\right|\leq\delta.$$
Let us consider $\al,\beta\in[0,1)$ such that
$\tau^+-\tau^--\beta=k\in\Z$, and
$U(\tau^+,0)-U(\tau^+-k,0)+\al\in\Z$. From \eqref{w3} we have
\begin{equation*}\begin{split}U(\tau^+,Y)&\leq U(\tau^+,0)+C_1\leq
U(\tau^+-k,Y)+2C_1+(U(\tau^+,0)-U(\tau^+-k,0))\\&\leq
U(\tau^+-k,Y)+2\lceil C_1
\rceil+(U(\tau^+,0)-U(\tau^+-k,0)+\al).\end{split}\end{equation*}
Since $\s(\cdot,y)$ and $W'(\cdot)$ are $\Z$-periodic, the
comparison principle for \eqref{wlip}  on the time interval
$[\tau^+,\tau^++T]$ implies that:
$$U(\tau^++T,Y)\leq
U(\tau^+-k+T,Y)+2\lceil C_1 \rceil+U(\tau^+,0)-U(\tau^+-k,0)+1 .$$
Choosing $Y=0$ in the previous inequality we get
\begin{equation*}\begin{split}U(\tau^++T,0)-U(\tau^+,0)&\leq
U(\tau^+-k+T,0)-U(\tau^+-k,0)+2\lceil C_1
\rceil+1\\&=U(\tau^-+\beta+T,0)-U(\tau^-+\beta,0)+2\lceil C_1
\rceil+1,\end{split}\end{equation*}and setting $t=\beta$ and
$\tau=\tau^-+T$ in \eqref{w4} and $\tau=\tau^-$ in \eqref{w5} we
finally obtain:
$$T\lam^+(T)\leq T\lam^-(T)+4\lceil C_1
\rceil+1+2\|W'\|_\infty+2\|\s\|_\infty+2\delta T.$$Since this is
true for any $\delta>0$, we conclude that:
$$|\lam^+(T)-\lam^-(T)|\leq \frac{4\lceil C_1
\rceil+1+2\|W'\|_\infty+2\|\s\|_\infty}{T}.$$ Now arguing as in
\cite{im} and \cite{imr}, we conclude that there exist
$\lim_{T\rightarrow+\infty}\lam^{\pm}(T)=:\lam$ and
$$|\lam^{\pm}(T)-\lam|\leq\frac{4\lceil C_1
\rceil+1+2\|W'\|_\infty+2\|\s\|_\infty}{T},$$ which implies that
$$|U(T,0)-\lam T|\leq 4\lceil C_1
\rceil+1+2\|W'\|_\infty+2\|\s\|_\infty,$$ and then, using
\eqref{w3} we get \eqref{w-lam}. The uniqueness of $\lam$ follows
from \eqref{w-lam}. Finally, \eqref{lambounds} is obtained from
\eqref{lam+-} as $T\rightarrow+\infty$.\finedim

 \subsection{Proof of Theorem \ref{ergodic}}
Let us consider the viscosity solution of \eqref{wlip} for
$\eta=0$. By Proposition \ref{ergodic2} we know that there exists
a unique $\lam$ such that $U(\tau,Y)/\tau$ converges to $\lam$ as
$\tau$ goes to $+\infty$ for any $Y\in \R^{N+1}$. Moreover, by
Proposition \ref{fromN+1toN},
 $U(\tau,y,0)$ is viscosity solution of \eqref{w}. Hence, the theorem follows immediately
 from the uniqueness of the viscosity solution of \eqref{w}.

\subsection{Proof of Proposition \ref{lipcorrect}}$\mbox{ }$\\
\noindent {\bf Step 1: Definition of $W^\pm_\eta$}\\
Let us denote  by $U^+_\eta$ the solution of \eqref{wlip} with
$a_0=C_1$, where $C_1$ is defined as in \eqref{w3}, and by
$U^-_\eta$ the solution of \eqref{wlip} with $a_0=0$. Let
$\lam_\eta^+=\lim_{\tau\rightarrow+\infty}\frac{U^+_\eta(\tau,Y)}{\tau}$
and
$\lam_\eta^-=\lim_{\tau\rightarrow+\infty}\frac{U^-_\eta(\tau,Y)}{\tau}$;
the existence of $\lam_\eta^+$ and $\lam_\eta^-$ is guaranteed by
Proposition \ref{ergodic2}. 
Now, we set 
$$W^+_\eta(\tau,Y):=U^+_\eta(\tau,Y)-\lam^+_\eta\tau$$
and
$$W^-_\eta(\tau,Y):=U^-_\eta(\tau,Y)-\lam^-_\eta\tau.$$
\noindent {\bf Step 2: Limits of $\lambda^\pm_\eta$}\\
By stability (see e.g. \cite{bi}), for
$\eta\rightarrow0^+$ the sequence $(U^+_\eta)_\eta$ converges to
$U$ solution of \eqref{wlip} with $\eta=0$. Moreover by
\eqref{lambounds} the sequence $(\lam^+_\eta)_\eta$ is bounded.
Take a subsequence $\eta_n\rightarrow0$ as $n\rightarrow+\infty$
such that $\lam^+_{\eta_n}\rightarrow\lam_\infty$ as
$n\rightarrow+\infty$. We want to show that $\lam_\infty=\lam$,
 where $\lam=\lim_{\tau\rightarrow+\infty}\frac{U(\tau,Y)}{\tau}$. By the proof of Theorem \ref{ergodic}, we know that $\lam$ is the same quantity
 defined in Theorem \ref{ergodic}. Using \eqref{w-lam}, we get
 \beqs\begin{split} |\lam-\lam_\infty|&\leq \left|\lam-\frac{U(\tau,0)}{\tau}\right|
 +\left|\frac{U(\tau,0)}{\tau}-\frac{U^+_{\eta_n}(\tau,0)}{\tau}\right|
 +\left|\frac{U^+_{\eta_n}(\tau,0)}{\tau}-\lam^+_{\eta_n}\right|+|\lam^+_{\eta_n}-\lam_\infty|
 \\& \leq \left|\lam-\frac{U(\tau,0)}{\tau}\right|+\left|\frac{U(\tau,0)}{\tau}-\frac{U^+_{\eta_n}(\tau,0)}{\tau}\right|
 +\frac{C_3}{\tau}+|\lam^+_{\eta_n}-\lam_\infty|\end{split}\eeqs
where $C_3$ does not depend on $n$. Then, passing to the limit
first as $n\rightarrow+\infty$ and then as
$\tau\rightarrow+\infty$, we obtain that $\lam=\lam_\infty$. This
implies that $\lam^+_\eta\rightarrow\lam$ as $\eta\rightarrow0$.

The same argument shows that $\lam^-_\eta\rightarrow\lam$ as
$\eta\rightarrow0$.\\
\noindent {\bf Step 3: $W^+_\eta$ and $W^-_\eta$ are respectively sub and supersolutions}\\
Since by \eqref{w3},
$C_0+\inf_{Y'}U^+_\eta(\tau,Y')-U^+_\eta(\tau,Y)\geq0$, $W^+_\eta$
is supersolution of \eqref{V} with $\lam=\lam^+_\eta$. Moreover,
by \eqref{w-lam}, $W^+_\eta$ is bounded on $\R^+\times\R^{N+1}$
uniformly w.r.t. $\eta$: $|W^+_\eta(\tau,Y)|\leq C_3$ for all
$(\tau,Y)\in\R^+\times\R^{N+1}$.\\
\noindent {\bf Step 4: regularity properties of $W^\pm_\eta$}\\
By \eqref{derivUyn+1}, $W^+_\eta$ is Lipschitz continuous w.r.t.
$y_{N+1}$ and $-1\leq \p_{y_{N+1}}W^+_\eta\leq
\frac{\|W''\|_\infty}{\eta}$. This implies that $W^+_\eta$ is also
a viscosity subsolution of

\begin{equation}\label{V^+sub}\left\{
  \begin{array}{ll}
    \lam^+_\eta+\p_{\tau} V=L+\I[V(\tau,\cdot,y_{N+1})]-W'(V+\lam^+_\eta\tau+P\cdot Y)+\s(\tau,y)\\
    \qquad\qquad\quad+C_1(\|W''\|_\infty+\eta)
    & \hbox{in } \R^+\times\R^{N+1}\\
    V(0,Y)=0 & \hbox{on }\R^{N+1}.
  \end{array}
\right.\end{equation}

By Proposition \ref{fromN+1toN}, $W^+_\eta$ is supersolution of
\eqref{V} and subsolution of \eqref{V^+sub} in $\R^+\times\R^N$
for any $y_{N+1}\in\R$.
 Then by Proposition \ref{regularityvisc}, $W^+_\eta$ is of class $C^\al$ w.r.t. $y$ uniformly in $y_{N+1}$ and $\eta$, for any $0<\al<1$.

 Similar arguments show that $W^-_\eta$ is subsolution of \eqref{V} with $\lam=\lam^-_\eta$, is bounded on $\R^+\times\R^{N+1}$,
  Lipschitz continuous w.r.t. $y_{N+1}$ with $-1\leq \p_{y_{N+1}}W^+_\eta\leq \frac{\|W''\|_\infty}{\eta}$ and H\"{o}lder continuous w.r.t. $y$.
  This concludes the proof of Proposition \ref{lipcorrect}.

 \subsection{Proof of Proposition \ref{Hprop}}
The continuity of $\overline{H}(p,L)$ follows from stability of
viscosity solutions of \eqref{w} (see e.g. \cite{bi}) and from
\eqref{w-lam}. Indeed, let $(p_n,L_n)$ be a sequence converging to
$(p_0,L_0)$ as $n\rightarrow+\infty$ and set
$\lam_n=\lam(p_n,L_n)$, $n\geq0$. By \eqref{w-lam}, we have for
any $\tau>0$
$$\left|\lam_n-\frac{w_n(\tau,y)}{\tau}\right|\leq
\frac{C_3}{\tau}.$$ Stability of viscosity solutions of \eqref{w}
implies that $w_n$ converges locally uniformly in $(\tau,y)$ to a
function $w_0$ which is a solution of \eqref{w} with
$(p,L)=(p_0,L_0)$. This implies that
$\limsup_{n\rightarrow+\infty}|\lam_n-\lam_0|\leq
\frac{2C_3}{\tau}$ for any $\tau>0$. Hence, we conclude that
$\lim_{n\rightarrow+\infty}\lam_n=\lam_0$.

Property (i) is an immediate consequence of \eqref{lambounds}.

The monotonicity in $L$ of $\overline{H}(p,L)$ comes from the
comparison principle.

Let us show (iii). Let $v$ be the solution of \eqref{v} and
$\lam=\lam(p,L)$. Set $\widetilde{v}(\tau,y):=v(\tau,-y)$. Remark
that $\I[\widetilde{v}(\tau,\cdot),y]=\I[v(\tau,\cdot),-y]$. If
$\s(\tau,\cdot)$ is even then $\widetilde{v}$ satisfies
\begin{equation*}
\begin{cases}
\lam+\p_{\tau}
\widetilde{v}=\I[\widetilde{v}(\tau,\cdot),y]+L-W'(\widetilde{v}+\lam
t-p\cdot y)+\s (\tau,y)&\text{in}\quad \R^+\times\R^N\\
\widetilde{v}(0,y)=0& \text{on}\quad \R^N.
\end{cases}
\end{equation*}By the uniqueness of $\lam$ we deduce that
$\lam(L,p)=\lam(L,-p)$, i.e. (iii).

Finally let us turn to (iv). Define
$\widetilde{v}(\tau,y):=-v(\tau,-y)$. If $W'(\cdot)$ and
$\s(\tau,\cdot)$ are odd functions, $\widetilde{v}$ satisfies
\begin{equation*}
\begin{cases}
-\lam+\p_{\tau}
\widetilde{v}=\I[\widetilde{v}(\tau,\cdot),y]-L-W'(\widetilde{v}-\lam
t+p\cdot y)+\s (\tau,y)&\text{in}\quad \R^+\times\R^N\\
\widetilde{v}(0,y)=0& \text{on}\quad \R^N.
\end{cases}
\end{equation*}As before, we conclude that
$\lam(-L,p)=-\lam(L,p)$, i.e. (iv).

\section{Smooth approximate correctors}\label{smoothcorsec}
In this section, we prove the existence of approximate correctors
that are smooth w.r.t. $y_{N+1}$, namely Proposition
\ref{apprcorrectors}. We first need the following lemma:
\begin{lem}\label{susum}
Let $u_1,u_2\in C_b(\R^+\times\R^N)$ be viscosity subsolutions
(resp., supersolutions) of \eqref{V} in $\R^+\times\R^N$, then
$u_1+u_2$ is viscosity subsolution (resp., supersolution) of
\begin{equation*}\left\{
  \begin{array}{ll}
    2\lam+\p_{\tau} v=2L+\I[v]-W'(u_1+P\cdot Y+\lam\tau)\\ \qquad\qquad-W'(u_2+P\cdot Y+\lam\tau)+2\s(\tau,y) & \hbox{in } \R^+\times\R^{N}\\
    v(0,y)=0 & \hbox{on }\R^{N}.
  \end{array}
\right.\end{equation*}
\end{lem}
For the proof  see Lemma 5.8 in \cite{cs}.

Next, let us consider a  positive smooth function
$\rho:\R\rightarrow\R$, with support in $B_1(0)$ and mass 1. We
define a sequence of mollifiers $(\rho_\delta)_\delta$ by
$\rho_\delta(s)=\frac{1}{\delta}\rho\left(\frac{s}{\delta}\right)$,
$s\in\R.$ Let $W^+_\eta$ (resp. $W^-_\eta$) be the Lipschitz
supersolution (resp. subsolution) of \eqref{V} with
$\lam=\lam^+_\eta$ (resp. $\lam=\lam^-_\eta$), whose existence  is
guaranteed by Proposition \ref{lipcorrect}. We define
\beq\label{V+-etadelta}
V^{\pm}_{\eta,\delta}(t,y,y_{N+1}):=W^{\pm}_\eta(t,y,\cdot)\star\rho_\delta(\cdot)
=\int_{\R}W^{\pm}_\eta(t,y,z)\rho_\delta(y_{N+1}-z)dz.\eeq
\begin{lem}\label{convollem}The functions $V^+_{\eta,\delta}$ and $V^-_{\eta,\delta}$ are respectively super and subsolution of
\begin{equation}\label{apprcorrequconv}\left\{
  \begin{array}{ll}
    \lam^{\pm}_\eta+\p_{\tau} V^{\pm}_{\eta,\delta}=L+\I[V^{\pm}_{\eta,\delta}(\tau,\cdot,y_{N+1})]+\s(\tau,y)
    \\ \qquad\qquad\quad\,\,\,\, -\int_{\R}W'(W^{\pm}_{\eta}(\tau,y,z)+p\cdot y+z+\lam^{\pm}_\eta\tau)\rho_\delta(y_{N+1}-z)dz
    & \hbox{in } \R^+\times\R^{N+1}\\
    V^{\pm}_\eta(0,Y)=0 & \hbox{on }\R^{N+1}.
  \end{array}
\right.\end{equation}
\end{lem}
\dim We prove the lemma for supersolutions. Let
$Q_h^e=e+[-h/2,h/2)$,
$\overline{\rho}_\delta(e,h)=\int_{Q_h^e}\rho_\delta(y)dy$
 and $$I_h(\tau,y,y_{N+1})=\sum_{e\in h\Z}W^+_\eta(\tau,y,y_{N+1}-e)\overline{\rho}_\delta(e,h).$$
 The function $I_h$ is a discretization of the convolution integral and by classical results, converges uniformly to
 $V^+_{\eta,\delta}$ as $h\rightarrow0$.
 By Proposition \ref{fromN+1toN}, $W^+_\eta$ is a viscosity supersolution of \eqref{V} also in $\R^+\times\R^N$. Then, by Lemma \ref{susum},
 for any $y_{N+1}\in\R$, $I_h(\tau,y,y_{N+1})$ is a supersolution of
 \begin{equation*}\left\{
  \begin{array}{ll}
    \lam^+_\eta+\p_{\tau} V=L+\I[V(\tau,\cdot,y_{N+1})]+ \s(\tau,y)\sum_{e\in h\Z}\overline{\rho}_\delta(e,h)
    \\ \qquad\qquad\quad-\sum_{e\in h\Z} W'(W^+_\eta(\tau,y,y_{N+1}-e)\\
    \qquad\qquad\quad+p\cdot y+(y_{N+1}-e)+\lam^+_\eta\tau)\overline{\rho}_\delta(e,h)
      & \hbox{in } \R^+\times\R^{N}\\
    V(0,y)=0 & \hbox{on }\R^{N}.
  \end{array}
\right.\end{equation*} Using the stability result for viscosity
solution of non-local equations, see \cite{bi}, we conclude that
$V^+_{\eta,\delta}$ is supersolution of \eqref{apprcorrequconv} in
$\R^+\times\R^N$ and hence also in $\R^+\times\R^{N+1}$.\finedim
\subsection{Proof of Proposition \ref{apprcorrectors}}
We first show that the functions $V^+_{\eta,\delta}$ and
$V^-_{\eta,\delta}$, defined in \eqref{V+-etadelta}, are
respectively super and subsolution of
\begin{equation}\label{equvetadelta}\left\{
  \begin{array}{ll}
    \lam^{\pm}_\eta+\p_{\tau} V^{\pm}_{\eta,\delta}=L+\I[V^{\pm}_{\eta,\delta}(\tau,\cdot,y_{N+1})]
    -W'(V^{\pm}_{\eta,\delta}+P\cdot Y+\lam^{\pm}_\eta\tau)\\ \qquad\qquad\quad+\s(\tau,y){\mp}C_{\eta,\delta}
    & \hbox{in } \R^+\times\R^{N+1}\\
    V^{\pm}_\eta(0,Y)=0 & \hbox{on }\R^{N+1},
  \end{array}
\right.\end{equation} where
$C_{\eta,\delta}=\|W''\|_\infty(2\delta\|W''\|_\infty/\eta+\delta)$.
Using \eqref{lipyn+1} and the properties of the mollifiers, we get
\beqs\begin{split}&\left|W'(V^{\pm}_{\eta,\delta}(\tau,y,y_{N+1})+p\cdot
y+y_{N+1}+\lam^{\pm}_\eta\tau)\right.
\\&\left.-\int_{\R}W'(W^{\pm}_{\eta}(\tau,y,z)+p\cdot
y+z+\lam^{\pm}_\eta\tau)\rho_\delta(y_{N+1}-z)dz\right|
\\&\leq\int_{\R}\left|W'(V^{\pm}_{\eta,\delta}(\tau,y,y_{N+1})+p\cdot y+y_{N+1}+\lam^{\pm}_\eta\tau)\right.\\&-
\left.W'(W^{\pm}_{\eta}(\tau,y,z)+p\cdot
y+z+\lam^{\pm}_\eta\tau)\right|\rho_\delta(y_{N+1}-z)dz
\\&\leq \|W''\|_\infty\int_{\R}\left[\left|V^{\pm}_{\eta,\delta}(\tau,y,y_{N+1})-W^{\pm}_{\eta}(\tau,y,z)\right|+|y_{N+1}-z|\right]
\rho_\delta(y_{N+1}-z)dz
\\&\leq\|W''\|_\infty\int_{\R}\left[\int_{\R}\left|W^{\pm}_{\eta}(\tau,y,r)-W^{\pm}_{\eta}(\tau,y,z)\right|\rho_\delta(y_{N+1}-r)dr
+|y_{N+1}-z|\right] \rho_\delta(y_{N+1}-z)dz
\\&\leq \|W''\|_\infty\int_{\R}\left[\int_{|y_{N+1}-r|\leq\delta}\frac{\|W''\|_\infty}{\eta}|r-z|\rho_\delta(y_{N+1}-r)dr+|y_{N+1}-z|\right]
\rho_\delta(y_{N+1}-z)dz
\\&\leq\|W''\|_\infty\int_{|y_{N+1}-z|\leq\delta}\left[\frac{\|W''\|_\infty}{\eta}(|y_{N+1}-z|+\delta)+|y_{N+1}-z|\right]
\rho_\delta(y_{N+1}-z)dz\\&
\leq\|W''\|_\infty\left(2\delta\frac{\|W''\|_\infty}{\eta}+\delta\right)\end{split}\eeqs
From this estimate and Lemma \ref{convollem}, we deduce that
$V^+_{\eta,\delta}$ and $V^-_{\eta,\delta}$ are respectively super
and subsolution of \eqref{equvetadelta}. Now, we choose
$\delta=\delta(\eta)$ such that
$\|W''\|_\infty(2\delta\|W''\|_\infty/\eta+\delta)=o_\eta(1)$ as
$\eta\rightarrow0$ and define \beqs
V^{\pm}_\eta(\tau,Y):=V^{\pm}_{\eta,\delta(\eta)}(\tau,Y).\eeqs
Then the functions $V^{\pm}_\eta$ are the desired super and
subcorrectors. Indeed, we have already shown that they are super
and subsolution
 of \eqref{apprcorrequ} with $\lam^+_\eta$ and $\lam^-_\eta$ satisfying \eqref{appcorr1}. Properties (i) and (ii) of Proposition
\ref{Hprop} can be shown as in the proof of the proposition.
Finally,  \eqref{appcorr2}, \eqref{appcorr3} and
\eqref{contrderivappcorr} easily follow from \eqref{lipcorrbound},
\eqref{lipyn+1},
 \eqref{contrderivappcorrlip} and the properties of the mollifiers.
 \finedim

\section{Appendix}\label{appendix}

\noindent {\bf Proof of Proposition \ref{regularityvisc}}\\
\noindent {\bf Heuristic arguments}\\
Before entering in the proof, let us start with an heuristic explanation.
Indeed, replacing $\partial_t u$ by $u$, we should get a similar result for a stationary solution of
$$\I[u] + g_2 \le u \le \I[u] + g_1$$
At a point $(x,y)$, with $x\neq y$, of supremum of
$$u(x)-u(y)-K|x-y|^\alpha$$
we have for $r>0$
$$\left\{\begin{array}{l}
u(x)\le g_1 +K\I^{1,r}[|\cdot-y|^\alpha,x]+\I^{2,r}[u,x]\\
\\
u(y)\ge g_2 -K\I^{1,r}[|x-\cdot|^\alpha,y]+\I^{2,r}[u,y]
\end{array}\right.$$

Setting $e=\frac{x-y}{|x-y|}$, $\varphi_\alpha(z)=|z|^\alpha$ and using the homogeneity of the functions, we get for $r=\sigma|x-y|$
$$\I^{1,r}[|\cdot-y|^\alpha,x]=-|x-y|^{\alpha-1}c^{\sigma}_\alpha=\I^{1,r}[|x-\cdot|^\alpha,y]\quad \mbox{with}\quad 
-c^{\sigma}_\alpha=\I^{1,\sigma}[\varphi_\alpha,e]$$

Therefore we get
$$u(x)-u(y)-K|x-y|^\alpha \le g_1-g_2 - K|x-y|^\alpha - 2K|x-y|^{\alpha-1}c^{\sigma}_\alpha+\I^{2,r}[u,x]-\I^{2,r}[u,y]$$
 By the maximal property of $(x,y)$, for  any $z\in \R^N$ we have  
$$u(x+z)-u(y+z)\leq u(x)-u(y)$$ which implies that 
$$\I^{2,r}[u,x]-\I^{2,r}[u,y]\leq 0$$ We conclude that 
$$u(x)-u(y)-K|x-y|^\alpha \le g_1-g_2 - K|x-y|^\alpha - 2K|x-y|^{\alpha-1}c^{\sigma}_\alpha$$
We can show that $c^{\sigma}_\alpha>0$, for $\sigma$ small enough and then an optimization on $|x-y|$ shows that for $K$ large enough, 
the right hand side is negative. This shows the H\"{o}lder estimate.

It turns out that the condition  $c^{\sigma}_\alpha>0$ is not satisfied for large values of $\sigma$.\\
\noindent {\bf Rigorous proof}\\
We use standard techniques from the theory of regularity of viscosity solutions of uniformly elliptic  second-order local operators, see \cite{il},  adapted to our context.

We argue by contradiction, assuming that  $u$ does not belong to
$C_x^{\al}(\R^+\times\R^N)$.  Let $u^{\ep,\ep'}$ and $u_{\ep,\ep'}$
be respectively the double-parameters sup and inf convolution of
$u$ in $\R^+\times\R^N$, i.e.
 $$u^{\ep,\ep'}(t,x)=\sup_{(s,y)\in \R^+\times\R^N}\left(u(s,y)-\frac{1}{2\ep}|x-y|^2-\frac{1}{2\ep'}(t-s)^2\right),$$
 $$u_{\ep,\ep'}(t,x)=\inf_{(s,y)\in \R^+\times\R^N}\left(u(s,y)+\frac{1}{2\ep}|x-y|^2+\frac{1}{2\ep'}(t-s)^2\right).$$
Then $u^{\ep,\ep'}$ is  semiconvex and is a subsolution of
$$\p_{t} u^{\ep,\ep'}=\I[u^{\ep,\ep'}(t,\cdot)]+g_1\quad\text{in}\quad (t_{\ep'},+\infty)\times\R^N$$
and $u_{\ep,\ep'}$ is semiconcave and is a supersolution of
$$\p_{t} u_{\ep,\ep'}=\I[u_{\ep,\ep'}(t,\cdot)]+g_2\quad\text{in}\quad (t_{\ep'},+\infty)\times\R^N,$$
where $t_{\ep'}\rightarrow0$ as $\ep'\rightarrow0$, see e.g.
Proposition III.2 in \cite{s}.
 
 Since $u$ is not H\"{o}lder continuous in $x$, there exists $\al\in(0,1)$ such that for any $K>0$ and $\ep,\ep'>0$
 
\beqs\begin{split}&\sup_{(t,x_1,x_2)\in\R^+\times \R^{2N}}u^{\ep,\ep'}(t,x_1)-u_{\ep,\ep'}(t,x_2)-K|x_1-x_2|^\al\\&
 \geq \sup_{(t,x_1,x_2)\in\R^+\times \R^{2N}}u(t,x_1)-u(t,x_2)-K|x_1-x_2|^\al\\&>0.\end{split}\eeqs

In order to make the supremum attained at some point, let us introduce  
 smooth positive functions
$\psi_1(t)$ and $\psi_2(x)$ with bounded first and second
derivatives  such that $\psi_1(t)\rightarrow+\infty$ as
$t\rightarrow+\infty$, $\psi_2(x)\rightarrow+\infty$ as
$|x|\rightarrow+\infty$ and there exists $K_0>0$ such that
$|\psi_2(x)|\leq K_0(1+\sqrt{|x|})$.  The last assumption on $\psi_2$ assures that $\I^2[\psi_2]$ is finite at any point. Then, for any $K>0$ and
$\ep$, $\ep'>0$ and $\beta>0$ small enough,
 the supremum on $\R^+\times \R^{2N}$ of the function 
 \beq\label{newlabellemmhold1}u^{\ep,\ep'}(t,x_1)-u_{\ep,\ep'}(t,x_2)-\phi(t,x_1,x_2),\eeq where
 $$\phi(t,x_1,x_2)=K|x_1-x_2|^\al+\beta\psi_1(t)+\beta\psi_2(x_1),$$ is positive and is attained at some point
 $(\ts,\xs_1,\xs_2)\in [0,+\infty)\times \R^{2N}$.  For  $\ep$, $\ep'$ small enough, $\xs_1\neq\xs_2$. Moreover, since   $u^{\ep,\ep'}(0,x)=u_{\ep,\ep'}(0,x)=0$ for any $x\in \R^N$, it turns out that actually $\ts>t_{\ep'}$.
  Remark that
\begin{equation}\label{eq::s11}
|\xs_1-\xs_2|\leq\left(\frac{2\sup_{(t,x)\in\R^+\times\R^N}|u(t,x)|}{K}\right)^\frac{1}{\al}.
\end{equation}

The function \eqref{newlabellemmhold1} is semiconvex, hence, by Aleksandrov's Theorem, twice differentiable almost eve\-ry\-where. 
Let us now introduce a perturbation of it, for which we can choose maximum points of twice differentiability. First
we transform $(\ts,\xs_1,\xs_2)$ into a strict
maximum point. In order to do that, we consider a smooth function
$h:\R^+\rightarrow\R$, with compact support, such that
$h(0)=0$ and $h(s)>0$ for $0<s<1$ and we set
$\theta(t,x_1,x_2)=h((t-\ts)^2)+h(|x_1-\xs_1|^2)+h(|x_2-\xs_2|^2)$.
Clearly $(\ts,\xs_1,\xs_2)$ is a strict maximum point of
$u^{\ep,\ep'}(t,x_1)-u_{\ep,\ep'}(t,x_2)-\phi(t,x_1,x_2)-\theta(t,x_1,x_2)$. Next we consider a smooth function $\chi:\R^N\rightarrow\R$ such
that $\chi(x)=1$ if $|x|\leq1/2$ and $\chi(x)=0$ for $|x|\geq1$.

By  Jensen's Lemma, see e.g. Lemma A.3 of \cite{cil}, for every small and positive $\delta$ there
exist $s^\delta\in\R,\,q_1^\delta,\,q_2^\delta\in\R^N$ with
$|s^\delta|,\,|q_1^\delta|,|q_2^\delta|\leq \delta$
such that the function
\begin{equation}\label{jensen}\Phi(t,x_1,x_2)=u^{\ep,\ep'}(t,x_1)-u_{\ep,\ep'}(t,x_2)-K|x_1-x_2|^\al-\varphi_1(t,x_1)-\varphi_2(x_2),\end{equation} where
$$\varphi_1(t,x_1)=\beta\psi_1(t)+\beta\psi_2(x_1)+h((t-\ts)^2)+h(|x_1-\xs_1|^2)+s^\delta t+\chi(x_1-\xs_1)q_1^\delta\cdot x_1,$$
$$\varphi_2(x_2)=h(|x_2-\xs_2|^2)+\chi(x_2-\xs_2)q_2^\delta\cdot x_2,$$
has a maximum at $(t^\delta,x_1^\delta,x_2^\delta)$, with
\begin{equation}\label{eq::s12}
|t_\delta-\ts|,\,|x_1^\delta-\xs_1|,\,|x_2^\delta-\xs_2|\leq\delta
\end{equation}

and $u^{\ep,\ep'}(t,x_1)-u_{\ep,\ep'}(t,x_2)$ is twice
differentiable at $(t^\delta,x_1^\delta,x_2^\delta)$. In
particular $u^{\ep,\ep'}$ is twice differentiable w.r.t. $x_1$ at
$(t^\delta,x_1^\delta)$ and $u_{\ep,\ep'}$ is twice differentiable
w.r.t. $x_2$ at $(t^\delta,x_2^\delta)$.  The function $\chi$ has been introduced to  make $\I^2[\varphi_1]$ and $\I^2[\varphi_2]$ finite. 

 For $\delta$ small
enough, we can assume $x_1^\delta\neq x_2^\delta$ and this will allow us to compute the derivatives of \eqref{jensen}. Since
$(t^\delta,x_1^\delta,x_2^\delta)$ is a maximum point, we have 
\begin{equation}\label{gradientzerohold}\begin{split}\nabla_{x_1}u^{\ep,\ep'}(t^\delta,x_1^\delta)=\nabla_{x_1}\varphi_1(t^\delta,x_1^\delta)+\al
K|x_1^\delta-x_2^\delta|^{\al-2}(x_1^\delta-x_2^\delta),\\
\nabla_{x_2}u_{\ep,\ep'}(t^\delta,x_2^\delta)=-\nabla_{x_2}\varphi_2(x_2^\delta)+\al
K|x_1^\delta-x_2^\delta|^{\al-2}(x_1^\delta-x_2^\delta).\end{split}\end{equation}
Moreover the inequalities
$$\Phi(t^\delta,x_1^\delta+z,x_2^\delta)\leq 
\Phi(t^\delta,x_1^\delta,x_2^\delta),$$
 
 $$\Phi(t^\delta,x_1^\delta,x_2^\delta+z)\leq 
\Phi(t^\delta,x_1^\delta,x_2^\delta),$$

$$\Phi(t^\delta,x_1^\delta+z,x_2^\delta+z)\leq 
\Phi(t^\delta,x_1^\delta,x_2^\delta),$$
for any $z\in\R^N$, with together \eqref{gradientzerohold}, give respectively: 

\beq\label{holdthmdis2}\begin{split}&u^{\ep,\ep'}(t^\delta,x_1^\delta+z)-u^{\ep,\ep'}(t^\delta,x_1^\delta)
-\nabla_{x_1}u^{\ep,\ep'}(t^\delta,x_1^\delta)\cdot z\\&\leq
\varphi_1(t^\delta,x_1^\delta+z)-\varphi_1(t^\delta,x_1^\delta)
-\nabla_{x_1}\varphi_1(t^\delta,x_1^\delta)\cdot z
\\&+K|x_1^\delta+z-x_2^\delta|^\al-K|x_1^\delta-x_2^\delta|^\al-\al
K|x_1^\delta-x_2^\delta|^{\al-2}(x_1^\delta-x_2^\delta)\cdot z,
\end{split}\end{equation}

\beq\label{holdthmdis3}\begin{split}&-(u_{\ep,\ep'}
(t^\delta,x_2^\delta+z)-u_{\ep,\ep'}(t^\delta,x_2^\delta)-\nabla_{x_2}u_{\ep,\ep'}(t^\delta,x_2^\delta)\cdot
z)\\&\leq\varphi_2(x_2^\delta+z)-\varphi_2(x_2^\delta)-\nabla_{x_2}\varphi_2(x_2^\delta)\cdot
z\\&+K|x_1^\delta-z-x_2^\delta|^\al-K|x_1^\delta-x_2^\delta|^\al+\al
K|x_1^\delta-x_2^\delta|^{\al-2}(x_1^\delta-x_2^\delta)\cdot z,
\end{split}\end{equation}
and for any  $r>0$
\begin{equation}\label{holdthmdis}\begin{split}&u^{\ep,\ep'}(t^\delta,x_1^\delta+z)-u^{\ep,\ep'}(t^\delta,x_1^\delta)
-\nabla_{x_1}u^{\ep,\ep'}(t^\delta,x_1^\delta)\cdot z
{\bf1}_{B_{r}}(z)\\&\leq u_{\ep,\ep'}
(t^\delta,x_2^\delta+z)-u_{\ep,\ep'}(t^\delta,x_2^\delta)-\nabla_{x_2}u_{\ep,\ep'}(t^\delta,x_2^\delta)\cdot
z
{\bf1}_{B_{r}}(z)\\&+\varphi_1(t^\delta,x_1^\delta+z)-\varphi_1(t^\delta,x_1^\delta)
-\nabla_{x_1}\varphi_1(t^\delta,x_1^\delta)\cdot z
{\bf1}_{B_{r}}(z)\\&
+\varphi_2(x_2^\delta+z)-\varphi_2(x_2^\delta)-\nabla_{x_2}\varphi_2(x_2^\delta)\cdot
z {\bf1}_{B_{r}}(z),\end{split}\end{equation} where
$B_{r}=B_{r}(0)$. The last inequality in particular implies that
\beq\label{holdthmdis4}\I^{2,r}[u^{\ep,\ep'}(t^\delta,\cdot),x_1^\delta]\leq
\I^{2,r}[u_{\ep,\ep'}(t^\delta,\cdot),x_2^\delta]
+\I^{2,r}[\varphi_1(t^\delta,\cdot),x_1^\delta]+\I^{2,r}[\varphi_2,x_2^\delta].\eeq

Next, in order to test, we need to double the time variables. Hence, for
$j>0$, let us consider the maximum point $(t^j,x_1^j,s^j,x_2^j)$
of the function
\beqs\begin{split}u^{\ep,\ep'}(t,x_1)-u_{\ep,\ep'}(s,x_2)-\Psi(t,x_1,x_2)-\frac{j}{2}|t-s|^2,\end{split}\eeqs
where \beqs
\Psi(t,x_1,x_2)=K|x_1-x_2|^\al+\varphi_1(t,x_1)+\varphi_2(x_2)
+|t-t^\delta|^2+|x_1-x_1^\delta|^2+|x_2-x_2^\delta|^2,\eeqs on
$Q_{\overline{\rho},\overline{\rho}}(t^\delta,x_1^\delta)\times
Q_{\overline{\rho},\overline{\rho}}(t^\delta,x_2^\delta)$, for
$\overline{\rho}>0$ sufficiently small. Standard arguments show
that
$(t^j,x_1^j,s^j,x_2^j)\rightarrow(t^\delta,x_1^\delta,t^\delta,x_2^\delta)$
as $j\rightarrow+\infty$. Hence for $j$ large enough there exists
$\rho>0$ such that $Q_{\rho,\rho}(t^j,x_1^j)\times
Q_{\rho,\rho}(s^j,x_2^j)\subset
Q_{\overline{\rho},\overline{\rho}}(t^\delta,x_1^\delta)\times
Q_{\overline{\rho},\overline{\rho}}(t^\delta,x_2^\delta)$  and
$x_1^j\neq x_2^j$. Testing, we get \beqs
j(t^j-s^j)+2(t^j-t^\delta)+\p_t\varphi_1(t^j,x_1^j)\leq
\I^{1,\rho}[\Psi(t^j,\cdot,x_2^j),x_1^j]+\I^{2,\rho}[u^{\ep,\ep'}(t^j,\cdot),x_1^j]+g_1,\eeqs
\beqs j(t^j-s^j)\geq
-\I^{1,\rho}[\Psi(t^j,x_1^j,\cdot),x_2^j]+\I^{2,\rho}[u_{\ep,\ep'}(s^j,\cdot),x_2^j]+g_2.\eeqs
Subtracting the two last inequalities, and then letting
$j\rightarrow+\infty$, we have \beqs\begin{split}
\p_t\varphi_1(t^\delta,x_1^\delta)&\leq
\I^{1,\rho}[\Psi(t^\delta,\cdot,x_2^\delta),x_1^\delta]+\I^{1,\rho}[\Psi(t^\delta,x_1^\delta,\cdot),x_2^\delta]
\\&+\I^{2,\rho}[u^{\ep,\ep'}(t^\delta,\cdot),x_1^\delta]-\I^{2,\rho}[u_{\ep,\ep'}(t^\delta,\cdot),x_2^\delta]+g_1-g_2.\end{split}\eeqs
Since $u^{\ep,\ep'}(t^\delta,\cdot)$ and
$u_{\ep,\ep'}(t^\delta,\cdot)$ are  twice differentiable
respectively at $x_1=x_1^\delta$ and $x_2=x_2^\delta$, we can pass
to the limit as $\rho\rightarrow0^+$ and obtain
 \beqs \p_t\varphi_1(t^\delta,x_1^\delta)
 \leq \I[u^{\ep,\ep'}(t^\delta,\cdot),x_1^\delta]-\I[u_{\ep,\ep'}(t^\delta,\cdot),x_2^\delta]+g_1-g_2.\eeqs Using
 \eqref{holdthmdis4}, we finally get
\begin{equation}\begin{split}\label{disjensen}\p_t\varphi_1(t^\delta,x_1^\delta)&\leq
\I^{1,r}[u^{\ep,\ep'}(t^\delta,\cdot),x_1^\delta]-\I^{1,r}[u_{\ep,\ep'}(t^\delta,\cdot),x_2^\delta]\\&
+\I^{2,r}[\varphi_1(t^\delta,\cdot),x_1^\delta]+\I^{2,r}[\varphi_2,x_2^\delta]+g_1-g_2.\end{split}\end{equation}

Next, let us estimate the term
$\I^{1,r}[u^{\ep,\ep'}(t^\delta,\cdot),x_1^\delta]-\I^{1,r}[u_{\ep,\ep'}(t^\delta,\cdot),x_2^\delta]$
and show that it contains a main negative part. For $0<\nu_0<1$,
let us denote
$$A_r:=\left\{z\in B_r(0)\,,\,|z\cdot (x_1^\delta-x_2^\delta)|\geq
\nu_0|z||x_1^\delta-x_2^\delta|\right\}.$$ Then
\beqs\begin{split}&\I^{1,r}[u^{\ep,\ep'}(t^\delta,\cdot),x_1^\delta]-\I^{1,r}[u_{\ep,\ep'}(t^\delta,\cdot),x_2^\delta]\\&=
\int_{A_r}[u^{\ep,\ep'}(t^\delta,x_1^\delta+z)-u^{\ep,\ep'}(t^\delta,x_1^\delta)
-\nabla_{x_1}u^{\ep,\ep'}(t^\delta,x_1^\delta)\cdot z \\&-(
u_{\ep,\ep'}
(t^\delta,x_2^\delta+z)-u_{\ep,\ep'}(t^\delta,x_2^\delta)-\nabla_{x_2}u_{\ep,\ep'}(t^\delta,x_2^\delta)\cdot
z )]\mu(dz)\\&+\int_{B_r\setminus
A_r}[...]\mu(dz)\\&=T_1+T_2.\end{split}\eeqs From \eqref{holdthmdis}
we have $$T_2\leq C.$$  Here and henceforth $C$ denotes various
positive constants independent of the pa\-ra\-me\-ters. Let us
estimate $T_1$. Using \eqref{holdthmdis2} and \eqref{holdthmdis3}, and successively making the change of variable $z\rightarrow -z$, we get the following estimate of $T_1$:
 \beqs\begin{split} T_1&\leq
\int_{A_r}[K|x_1^\delta+z-x_2^\delta|^\al-K|x_1^\delta-x_2^\delta|^\al-\al
K|x_1^\delta-x_2^\delta|^{\al-2}(x_1^\delta-x_2^\delta)\cdot
z]\mu(dz)+C\\&+\int_{A_r}[K|x_1^\delta-z-x_2^\delta|^\al-K|x_1^\delta-x_2^\delta|^\al+\al
K|x_1^\delta-x_2^\delta|^{\al-2}(x_1^\delta-x_2^\delta)\cdot
z]\mu(dz)\\&=2\int_{A_r}[K|x_1^\delta+z-x_2^\delta|^\al-K|x_1^\delta-x_2^\delta|^\al-\al
K|x_1^\delta-x_2^\delta|^{\al-2}(x_1^\delta-x_2^\delta)\cdot
z]\mu(dz)+C\\& \leq\al K\int_{A_r}\sup_{|t|\leq
1}\{|x_1^\delta-x_2^\delta+tz|^{\al-4}(|x_1^\delta-x_2^\delta+tz|^2|z|^2\\&-(2-\al)[(x_1^\delta-x_2^\delta+tz)\cdot
z]^2)\}\mu(dz)+C.\end{split}\eeqs Let us fix
$r=\sigma|x_1^\delta-x_2^\delta|$, $\sigma>0$, then for $z\in A_r$
$$|x_1^\delta-x_2^\delta+tz|\leq
(1+\sigma)|x_1^\delta-x_2^\delta|,$$
$$|(x_1^\delta-x_2^\delta+tz)\cdot
z|\geq|(x_1^\delta-x_2^\delta)\cdot z|-|z|^2\geq
\left(\nu_0-\sigma\right)|x_1^\delta-x_2^\delta||z|.$$ Let us
choose $0<\sigma<\nu_0<1$ such  that
$$C_0:=-(1+\sigma)^2+(2-\al)(\nu_0-\sigma)^2>0,$$then \beqs T_1\leq
-CC_0 K|x_1^\delta-x_2^\delta|^{\al-2}
\int_{A_r}|z|^2\mu(dz)+C.\eeqs By homogeneity
$$\int_{A_r}|z|^2\mu(dz)=Cr.$$
Then, we conclude $$T_1\leq
-CC_0K|x_1^\delta-x_2^\delta|^{\al-2}r+C\leq
-CC_0K|x_1^\delta-x_2^\delta|^{\al-1}+C,$$
and   from \eqref{disjensen}  
              \begin{equation*}\begin{split}C C_0K|x_1^\delta-x_2^\delta|^{\al-1}&\leq
              -\p_t\varphi_1(t^\delta,x_1^\delta)+g_1-g_2+C\\&
              +\I^{2,r}[\varphi_1(t^\delta,\cdot),x_1^\delta]+\I^{2,r}[\varphi_2,x_2^\delta]\\&\leq
g_1-g_2+C.\end{split}\end{equation*} Letting $\delta$ go to 0, from the previous inequalities and \eqref{eq::s12} we finally obtain
 \begin{equation*}K|\xs_1-\xs_2|^{\al-1}\leq C,\end{equation*} where $C$ is independent of $K$.
              This is a contradiction for  $K$ large enough, because of \eqref{eq::s11}, hence $u\in C^\al_x(\R^+\times\R^N)$.
\finedim

\noindent {\bf Proof of Proposition \ref{pro::s19}}\\ 
Let us define the functions
$V_1(\tau,Y):=e^{-k\tau}U_1(\tau,Y)$ and
$V_2(\tau,Y):=e^{-k\tau}U_2(\tau,Y)$, where $k:=\|W''\|_\infty+1$.
It is easy to see that $V_1$ and $V_2$ are respectively sub and
supersolution of
\begin{equation}\label{compeq}\left\{
  \begin{array}{ll}
    \p_{\tau} V=Le^{-k\tau}+\I[V(\tau,\cdot,y_{N+1})]+g(\tau,Y,V)\\
    \quad\quad\,+\eta[a_0+e^{k\tau}(\inf_{Y'}V(\tau,Y')-V(\tau,Y))]|\p_{y_{N+1}}V+e^{-k\tau}|
    & \hbox{in } \R^+\times\R^{N+1}\\
    V(0,Y)=0 & \hbox{on }\R^{N+1},
  \end{array}
\right.\end{equation}

where $g(\tau,Y,V)=-e^{-k\tau}W'(e^{k\tau}V+P\cdot
Y)-kV+e^{-k\tau}\s(\tau,y)$. Remark that, by the choice of $k$,
\beq\label{compar1}g(\tau,Y,V_1)-g(\tau,Z,V_2)\leq
-(V_1-V_2)+e^{-k\tau}(\|W''\|_\infty|P|+\|\s'\|_\infty)|Y-Z|.\eeq
To prove the comparison between $U_1$ and $U_2$, it suffices to
show that  $V_1(\tau,Y)\leq V_2(\tau,Y)$ for all $(\tau,Y)\in
(0,T)\times \R^{N+1}$ and for any $T>0$.

Suppose by contradiction that
$M=\sup_{(\tau,Y)\in(0,T)\times\R^{N+1}}(V_1(\tau,Y)-V_2(\tau,Y))>0$
for some $T>0$. Define for small $\nu_1,\nu_2,\beta,\delta>0$ the
function $\phi\in C^2((\R^+\times\R^{N+1})^2)$ by
$$\phi(\tau,Y,s,Z)=\frac{1}{2\nu_1}|\tau-s|^2+\frac{1}{2\nu_2}|Y-Z|^2+\beta\psi(Y)+\frac{\delta}{T-\tau},$$ where $\psi$ is defined as the function
$\psi_2$ in the proof of Proposition \ref{regularityvisc}. The
supremum of $V_1(\tau,Y)-V_2(s,Z)-\phi(\tau,Y,s,Z)$ is attained at
some point
 $(\tas,\Ys,\os,\Zs)\in ((0,T)\times\R^{N+1})^2$. Standard arguments show
 that, because $U_1$ and $U_2$ are assumed bounded
 \beqs (\tas,\Ys,\os,\Zs)\rightarrow(\widehat{\tau},\widehat{\tau},\widehat{Y},\widehat{Z})\quad\text{as }\nu_1\rightarrow0,\eeqs
\beqs V_1(\tas,\Ys)\rightarrow V_1(\widehat{\tau},\widehat{Y}),\,
V_2(\os,\Zs)\rightarrow
V_2(\widehat{\tau},\widehat{Z})\quad\text{as
}\nu_1\rightarrow0,\eeqs
 where
 $(\widehat{\tau},\widehat{Y},\widehat{Z})$ is a maximum point of $V_1(\tau,Y)-V_2(\tau,Z)-\frac{1}{2\nu_2}|Y-Z|^2-\beta\psi(Y)-\frac{\eta}{T-\tau}$.
 Moreover, it is easy to see that
 \beqs \limsup_{\nu_1\rightarrow0}\inf_{Y'}V_1(\tas,Y')\leq \inf_{Y'}V_1(\widehat{\tau},Y'),\,
 \liminf_{\nu_1\rightarrow0}\inf_{Y'}V_2(\os,Y')\geq\inf_{Y'}V_2(\widehat{\tau},Y').\eeqs

Since $V_1$ and $V_2$ are respectively sub and supersolution of
\eqref{compeq}, for any $r>0$ we have
\beq\label{comp4}\begin{split}&\frac{\delta}{(T-\tas)^2}+\frac{\tas-\os}{\nu_1}\\&\leq
Le^{-k\tas}+\frac{C_Nr}{\nu_2}
+\beta\I^{1,r}[\psi(\cdot,\ys_{N+1}),\ys]+\I^{2,r}[V_1(\tas,\cdot,\ys_{N+1}),\ys]+g(\tas,\Ys,V_1(\tas,\Ys))\\&
+\eta[a_0+e^{k\tas}(\inf_{Y'}V_1(\tas,Y')-V_1(\tas,\Ys))]\left|\frac{\ys_{N+1}-\zs_{N+1}}{\nu_2}+\beta
\p_{y_{N+1}}\psi(\Ys)+ e^{-k\tas}\right|\end{split}\eeq and
\beq\label{comp5}\begin{split}\frac{\tas-\os}{\nu_1}&\geq
Le^{-k\os}-\frac{C_Nr}{\nu_2}
+\I^{2,r}[V_2(\os,\cdot,\zs_{N+1}),\zs]+g(\os,\Zs,V_2(\os,\Zs))\\&
+\eta[a_0+e^{k\os}(\inf_{Y'}V_2(\os,Y')-V_2(\os,\Zs))]\left|\frac{\ys_{N+1}-\zs_{N+1}}{\nu_2}+e^{-k\os}\right|,\end{split}\eeq
where $C_N$ is a constant depending on the dimension $N$. Since
$(\tas,\Ys,\os,\Zs)$ is a maximum point, we have
\begin{equation*}\begin{split}
&V_1(\tas,\ys+x,\ys_{N+1})-V_1(\tas,\Ys)\leq
V_2(\overline{s},\zs+x,\zs_{N+1})-V_2(\overline{s},\Zs)+\beta[\psi(\ys+x,\ys_{N+1})-\psi(\Ys)],\end{split}\end{equation*}
for any $x\in\R^N$, which implies that for any $r>0$
\begin{equation*}\I^{2,r}[V_1(\tas,\cdot,\ys_{N+1}),\ys]\leq \I^{2,r}[V_2(\os,\cdot,\zs_{N+1}),\zs]
+\beta\I^{2,r}[\psi(\cdot,\ys_{N+1}),\ys].\end{equation*} Then,
subtracting \eqref{comp4} with \eqref{comp5} and letting
$r\rightarrow0^+$, we get \beqs\begin{split}
\frac{\delta}{(T-\tas)^2}&\leq
L(e^{-k\tas}-e^{-k\os})+\beta\I[\psi(\cdot,\ys_{N+1}),\ys]
+g(\tas,\Ys,V_1(\tas,\Ys))-g(\os,\Zs,V_2(\os,\Zs))\\&
+\eta[a_0+e^{k\tas}(\inf_{Y'}V_1(\tas,Y')-V_1(\tas,\Ys))]\left|\frac{\ys_{N+1}-\zs_{N+1}}{\nu_2}+\beta
\p_{y_{N+1}}\psi(\Ys)+e^{-k\tas}\right|\\&
-\eta[a_0+e^{k\os}(\inf_{Y'}V_2(\os,Y')-V_2(\os,\Zs))]\left|\frac{\ys_{N+1}-\zs_{N+1}}{\nu_2}+e^{-k\os}\right|.
\end{split}\end{equation*}
Next, letting $\nu_1\rightarrow0$ and using \eqref{compar1}, we
obtain \beq\label{comp6}\begin{split}
&\frac{\delta}{(T-\widehat{\tau})^2}\\&\leq
\beta\I[\psi(\cdot,\widehat{y}_{N+1}),\widehat{y}]
-(V_1(\widehat{\tau},\widehat{Y})-V_2(\widehat{\tau},\widehat{Z}))+e^{-k\widehat{\tau}}(\|W''\|_\infty|P|
+\|\s'\|_\infty)|\widehat{Y}-\widehat{Z}|+C\beta\\& +\eta
e^{k\widehat{\tau}}[\inf_{Y'}V_1(\widehat{\tau},Y')-\inf_{Y'}V_2(\widehat{\tau},Y')-(V_1(\widehat{\tau},\widehat{Y})-V_2(\widehat{\tau},\widehat{Z}))]
\left|\frac{\widehat{y}_{N+1}-\widehat{z}_{N+1}}{\nu_2}+e^{-k\widehat{\tau}}\right|.\end{split}\end{equation}
It is easy to prove that \beq\label{comp3}
\liminf_{(\beta,\delta)\rightarrow(0,0)}(V_1(\widehat{\tau},\widehat{Y})-V_2(\widehat{\tau},\widehat{Z}))\geq
M\eeq and \beqs \frac{|\widehat{Y}-\widehat{Z}|^2}{\nu_2}\leq
C,\eeqs where $C$ is independent of $\beta$ and $\delta$. Up to
subsequence, $\widehat{\tau}\rightarrow \tau_0\in[0,T]$ as
$(\beta,\delta)\rightarrow(0,0)$ and by \eqref{comp3}, we have
\beqs\begin{split}&
\limsup_{(\beta,\delta)\rightarrow(0,0)}[\inf_{Y'}V_1(\widehat{\tau},Y')-\inf_{Y'}V_2(\widehat{\tau},Y')
-(V_1(\widehat{\tau},\widehat{Y})-V_2(\widehat{\tau},\widehat{Z}))]
\\&\leq \inf_{Y'}V(\tau_0,Y')-\inf_{Y'}V_2(\tau_0,Y')
-\sup_{Y'}(V_1(\tau_0,Y')-V_2(\tau_0,Y'))\\&\leq0.\end{split}\eeqs

Then, passing to the limit first as
$(\beta,\delta)\rightarrow(0,0)$ and then as $\nu_2\rightarrow0$
in \eqref{comp6} we finally get the contradiction:
$$M\leq0,$$and this concludes the proof of the comparison theorem. \finedim


\end{document}